%% file: adaptive_multi_index_collocation_rev1.tex
\pdfoutput=1
\def\SNL{Optimization and Uncertainty Quantification, Sandia National Laboratories, Albuquerque, NM, 87123}
\def\MICH{Department of Aerospace Engineering  University of Michigan  Ann Arbor, MI, 48109}

\ifx \journal\undefined
\documentclass[10pt]{uq-report}
\usepackage{authblk}

\author[1]{J.D. Jakeman}
\affil[1]{\SNL}
\author[1]{M.S. Eldred}
\author[1]{G. Geraci}
\author[2]{A. Gorodetsky}
\affil[2]{\MICH}
\corauthor{J.D. jakeman}
\coremail{jdjakem@sandia.gov}
\shortauthor{J.D Jakeman, M.Eldred, G. Geraci, A. Gorodetsky}
\shorttitle{Adaptive MISC for UQ and SA}
\funding{See acknowledgements}
\reportname{Preprint}
\else
\documentclass[num-refs]{wiley-article}
\newtheorem{assumption}{Assumption}

\papertype{Original Article}
\author[1\authfn{1}]{J.D. Jakeman PhD}
\author[1\authfn{2}]{M.S. Eldred PhD}
\author[1\authfn{2}]{G. Geraci PhD}
\author[2\authfn{2}]{A. Gorodetsky PhD}
\affil[1]{Optimization and Uncertainty Quantification, Sandia National Laboratories, Albuquerque, NM, 87123, USA}
\affil[2]{Department of Aerospace Engineering  University of Michigan  Ann Arbor, MI, 48109, USA}
\corraddress{J.D. Jakeman PhD, Optimization and Uncertainty Quantification, Sandia National Laboratories, Albuquerque, NM, 87123, USA}
\corremail{jdjakem@sandia.gov}
\fundinginfo{All authors were supported by DARPA EQUiPS.}
\runningauthor{J.D. Jakeman et al.}
\fi

\title{Adaptive Multi-index Collocation for \rev{Uncertainty Quantification} and Sensitivity Analysis}

\include{defs}
\usepackage[table]{xcolor}
\usepackage[export]{adjustbox}
\usepackage{graphbox}
\usepackage{slashbox}
\usepackage{floatrow}
\newfloatcommand{capbtabbox}{table}[][\FBwidth]
\usepackage{algorithm}
\usepackage{algorithmicx,algpseudocode}
\usepackage[textsize=tiny]{todonotes}
\usepackage{wrapfig}
\usepackage{mathtools}

\usepackage{amsmath,amsfonts,amssymb,graphicx,url,hyperref}

\newcolumntype{H}{>{\setbox0=\hbox\bgroup}c<{\egroup}@{}}

\begin{document}

\maketitle

\begin{abstract}
  In this paper, we present an adaptive algorithm to construct response surface approximations of high-fidelity models using a hierarchy of lower fidelity models. Our algorithm is based on multi-index stochastic collocation and automatically balances physical discretization error and response surface error to construct an approximation of model outputs. This surrogate can be used for uncertainty quantification (UQ) and sensitivity analysis (SA) at a fraction of the cost of a purely high-fidelity approach.  We demonstrate the effectiveness of our algorithm on a canonical test problem from the UQ literature and a complex multi-physics model that simulates the performance of an integrated nozzle for an unmanned aerospace vehicle. We find that, \rev{when the input-output response is sufficiently smooth}, our algorithm produces approximations that can be \rev{over two} orders of magnitude more accurate than single fidelity approximations for a fixed computational budget.
  \ifdefined\journal
  \fi
\end{abstract}
\section{Introduction}
\label{sec-1}
Quantifying uncertainty in simulation-based prediction of engineered systems is essential for credible certification and design. Traditional sampling-based approaches \rev{for quantifying uncertainty}, such as Monte Carlo (MC), require repeated evaluations at varying parameter realizations. \rev{Consequently, MC methods} are often computationally intractable for high-fidelity simulation models because they require a prohibitively large number of evaluations to obtain moderate accuracy in response statistics. \rev{In contrast,} surrogate methods such as polynomial chaos \cite{Ghanem_S_book_1991,Sudret_RESS_2008,Xiu_K_SISC_2002,Arnst_GPR_IJNME_2013}, Gaussian processes \cite{Rasmussen_EW_book_2005,Sacks_WMW_SS_1989}, low-rank decompositions \cite{Doostan_VI_CMAME_2013,Gorodetsky_J_JCP_2018,Oseledets2011,Gorodetsky2015}, and sparse grid interpolation \cite{Nobile_TW_SINUM_08,Xiu_H_SISC_2005,Jakeman_R_SGA_2013,Agarwal_A_IJNME_2010} can be used to build an approximation of the input-output response, often at a fraction of the cost of MC sampling. Once the surrogate has been constructed, various uncertainty quantification (UQ) tasks, such as sensitivity analysis, density estimation, etc., can then be performed on the approximation at negligible cost.\footnote{\rev{Reduced order models can also be used to construct surrogates and have been applied successfully for UQ on many applications \cite{Cui_MW_IJNME_2014,Soize_F_IJNME_2017,Carlberg_IJNME_2014}. These methods do not construct response surface approximations, but rather solve the governing equations on a reduced basis.}}

Despite the improved efficiency of surrogate methods relative to MC sampling, building a surrogate can still be prohibitively expensive for high-fidelity simulation models. Fortunately, a selection of models of varying fidelity and computational cost are typically available for \rev{many applications. For example, aerospace models span} fluid dynamics, structural and thermal response, control systems, etc. 
Leveraging an ensemble of models can facilitate significant reductions in the overall computational cost of UQ, by integrating the predictions of quantities of interest (QoI) from multiple sources.


%

Multi-fidelity methods utilize an ensemble of models, enriching a small number of high-fidelity simulations with larger numbers of simulations from models of varying prediction accuracy and reduced cost, to enable greater exploration and resolution of uncertainty while maintaining deterministic prediction accuracy. The effectiveness of multi-fidelity approaches depends on the ability to identify and exploit relationships among models within the ensemble.
The relationships among models within a model ensemble vary greatly, and most existing approaches focus on exploiting a specific type of structure for a presumed model sequence. For example, \cite{Perdikaris_PRSLA_2016,Kennedy_O_B_2000,LeGratiet_G_IJUQ_2014,Narayan_GX_SISC_2013} build approximations that leverage a hierarchy of models of increasing fidelity, with varying physics and/or numerical discretizations. Multi-level \cite{Teckentrup_JWG_SIAMUQ_2015} and multi-index \cite{HajiAli_NTT_CMAME_2016} also leverage a set of models of increasing fidelity, with the additional assumption that a model sequence forms a convergent hierarchy. Such a hierarchy can be formed from a sequence of finite element discretizations to the solutions of a partial differential equation (PDE)\rev{, for example}, which converges to the exact solution of the governing equations as the finite element mesh is refined. Some attention has been given to building multi-fidelity approximations using models that do not admit a strict ordering of fidelity \rev{\cite{Gorodetsky_GEJ_2018,Peherstorfer_WG_SIAM_2016}}; however literature in this area is limited.

In this paper, we focus on multi-level/multi-index surrogate methods, observing that access to a sequence of numerical model discretizations can be more practical than access to an ensemble of models of varying physics fidelity. \rev{Multi-level methods}~\cite{Ng_E_AIAA_2012,Eldred_NBD_UQHandbook_2016,Teckentrup_JWG_SIAMUQ_2015} build the multi-fidelity surrogate as a linear combination of approximations of the differences between the QoI of models of increasing fidelity/discretization. \rev{Resources are allocated to each model level, controlled by a single hyper-parameter, specifying the level of discretization for example, in a manner that balances computational cost with increasing accuracy.}
In many applications, however, multiple hyper-parameters may control the model discretization, such as the mesh and time step sizes. In these situations, it may not be clear how to construct a one-dimensional hierarchy represented by a scalar hyper-parameter. To overcome this limitation, a generalization of multi-level collocation, referred to as multi-index stochastic collocation (MISC), was developed to deal with multivariate hierarchies with multiple refinement hyper-parameters \cite{HajiAli_NTT_CMAME_2016}.

In this paper we will present an adaptive extension of MISC and use this algorithm to reduce the cost of quantifying uncertainty in the performance of a supersonic jet engine nozzle. The main contributions of this paper are:
\begin{itemize}
\item The extension of MISC to function approximation. The original MISC papers focused on estimation of moments using quadrature.
\item The development of an adaptive multi-index collocation approach based upon sparse grid approximation. The seminal multi-index collocation paper used {\it a priori} estimates of model and parameter importance to allocate samples across the model ensemble. Our algorithm adaptively balances \rev{discretization error introduced by numerically solving the governing equations}, response surface error, and parameter importance.
\item Formulation of a multi-index method for variance based sensitivity analysis.
\item Application of multi-index collocation to an engineering application of practical importance, specifically UQ of a nozzle of an unmanned aerospace vehicle.
\end{itemize}



\rev{Multi-level/multi-index surrogate methods are closely related to many multi-level/multi-fidelity sampling algorithms~\cite{Cliffe_GST_CVS_2011,Geraci_EI_AIAA_2017,Giles_OR_2008,Haji_NT_NM_2016,Ng_Thesis_MIT_2013, Peherstorfer_WG_SIAM_2016,Gorodetsky_GEJ_2018}. These sampling algorithms leverage correlation between the outputs of multiple models to reduce the variance in statistical estimators of quantities such as expectation. This variance reduction can result in orders of magnitude reduction in the computational cost of quantifying uncertainty, but like traditional MC sampling, the error decreases slowly as the number of samples is increased which can still render this approach infeasible in some contexts. As a rough guide, multi-fidelity surrogate methods are best applied to models which with moderate dimensionality and smooth input-output responses and sampling based methods are better suited to models with large numbers of parameters and for QoI with low regularity.}


The remainder of this paper is organized as follows. Section \ref{sec:models} introduces two motivating examples: an advection-diffusion model problem and an aero-thermo-structural analysis of the performance of a supersonic jet engine nozzle. Section \ref{sec:discretization} presents the numerical discretization used to solve the aforementioned models, and introduces the sparse grid interpolation method which is the foundation of our approach.  The adaptive multi-index collocation method is presented in Section \ref{sec:misc}, and Section \ref{sec:sensitivity-analysis} describes how to use MISC for global sensitivity analysis. Finally, in section \ref{sec:results}, we highlight the strengths of the adaptive multi-index method (AMISC) as demonstrated on the advection-diffusion model problem and the engineering-scale model of the supersonic nozzle.

\section{Problem setup}
\label{sec-2}
\label{sec:models}
We seek \rev{to quantify uncertainty} in a broad class of stochastic partial differential equations (PDE).
Let $\physdom\subset\mathbb{R}^\ell$, for $\ell=1,2,$ or $3$ \rev{and $\text{dim}(\physdom)=\ell$}, be a physical domain; $T>0$ be a real number; and $\rvvsupp\rev{\subseteq}
\mathbb{R}^{\nv}$, for $\nv\geq 1$ \rev{and $\text{dim}(\rvvsupp)=\nv$}, be the stochastic space. The PDE is defined as 
\begin{align} \label{eq:govern}
\left\{
\begin{array} {ll}
u_t(\dvv,t,\rVv) = \mathcal{L}(u), &
\\
\mathcal{B}(u(\dvv, t, \rVv)\rev{)} = 0, & 
 \\
u(\dvv, 0, \rVv) = u_0(\dvv, \rVv), & 
\end{array}
\right.
\end{align}
where $t$ is the time, $\dvv=(\dv_1,\dots,\dv_\ell)$ is the physical coordinate, $\mathcal{L}$ is a (nonlinear) differential operator, $\mathcal{B}$ is the boundary condition operator, $u_0$ is the initial
condition, and $\rVv=(\rV_1,\dots,\rV_{\nv})\in \rvvsupp$  are a set of random variables characterizing the random inputs to the governing equation.
The solution $u \in V$ is a vector-valued stochastic quantity
\begin{equation} \label{u}
u:\bar{\physdom}\times[0,T]\times \rvvsupp\to\mathbb{R}^{n},
\end{equation}
\rev{for some suitable function space $V$ and} where $\bar{\physdom} = \physdom \times \partial \physdom$ is the closure of the interior domain. As an example, $u(\dvv,t, \rvv)$ can be a vector of temperatures, pressures, and velocities for a specific location in the domain $\dvv$, specific time $t$, and for a specific realization $\rvv$ of the stochastic variable $\rVv$.

Often, one is more interested in quantifying uncertainty in a particular functional of the solution, called the quantity of interest (QoI), rather than the full solution. Let 
\begin{align}
\label{eq:qoi-functional}
F\left[u\right]: V \rightarrow \mathbb{R}^{\nq}, && \rev{\nq>0}.
\end{align}
be such a function, \rev{where $F$ is typically a continuous and bounded functional}. Then, we are interested in estimating statistics of the function
\begin{align}\label{eq:qoi-function}
  f(\rvv)=F[u(\cdot,\rvv)] : \Gamma \to \mathbb{R}^{\nq}.
\end{align}
which assigns the value of the quantity of interest $F[u]$ to each realization of the random variables $\rvv$.
Expression ~\eqref{eq:qoi-function} allows us to express the QoI as a function of the random variables without needed to explicitly state the dependence of the QoI on the PDE solution.



\subsection{Advection diffusion equation}
\label{sec-2-1}
\rev{As a stochastic model problem, we will use the following transient advection-diffusion PDE, in two spatial dimensions ($\physdom=[0,1]^2$), to highlight the strengths and weaknesses of our proposed approach,
\begin{align}
\label{eq:advection-diffusion}
\frac{\partial u}{\partial t}(\dvv,t,\rvv) + \nabla u(\dvv,t,\rvv)-\nabla\cdot\left[k(\dvv,\rvv) \nabla u(\dvv,t,\rvv)\right] = g(\dvv,t)& &
(\dvv,t,\rvv)\in\physdom\times [0,1]\times\rvvsupp \\\nonumber
u(\dvv,t,\rvv)=0 & & (\dvv,t,\rvv)\in \partial \physdom\times[0,1]\times\rvvsupp \\\nonumber
\end{align}
with forcing $g(\dvv,t)=(1.5+\cos(2\pi t))\cos(\dv_1)$, and subject to the initial condition $u(\dvv,0,\rvv)=0$. Following \cite{Nobile_TW_SIAMNA_2008}, we model the diffusivity $k$ as a random field represented by the
Karhunen-Lo\'{e}ve (like) expansion (KLE)
\begin{align}
\label{eq:diffusivityZ}
\log(k(\dvv,\rvv)-0.5)=1+\rv_1\left(\frac{\sqrt{\pi L}}{2}\right)^{1/2}+\sum_{k=2}^d \lambda_k\phi(\dvv)\rv_k,
\end{align}
with
\begin{align}
  \lambda_k=\left(\sqrt{\pi L}\right)^{1/2}\exp\left(-\frac{(\lfloor\frac{k}{2}\rfloor\pi L)^2}{4}\right) k>1,  & &  \phi(\dvv)=
    \begin{cases}
      \sin\left(\frac{(\lfloor\frac{k}{2}\rfloor\pi \dv_1)}{L_p}\right) & k \text{ even}\,,\\
      \cos\left(\frac{(\lfloor\frac{k}{2}\rfloor\pi \dv_1)}{L_p}\right) & k \text{ odd}\,.
    \end{cases}
\end{align}
where $L_p=\max(1,2L_c)$, $L=\frac{L_c}{L_p}$ and $L_c=0.5$.

In this example we truncate \eqref{eq:diffusivityZ} at $\nv=10$ terms and choose $\rV_k$, $k\in\{1\ldots,\nv\}$ to be uncorrelated independent and identical uniform random variables distributed on $[-\sqrt{3},\sqrt{3}]$. Using 10 terms captures $99.99999996$ percent of the total energy of the untruncated expansion.
We choose a random field which is effectively one-dimensional so that the error in the finite element solution is more sensitive to refinement of the mesh in the $\dv_1$-direction than to refinement in the $\dv_2$-direction.

The advection diffusion equation \eqref{eq:advection-diffusion} is solved using linear finite elements and implicit backward-Euler timestepping. In the following we will show how solving the PDE with varying numbers of finite elements and timesteps can reduce the cost of approximating the quantity of interest
\begin{align}\label{eq:advec-diff-qoi}f(\rvv)=\int_\physdom u(\rvv)\frac{1}{2\pi\sigma^2}\exp\left(-\frac{\lVert\dvv-\dvv^\star \rVert_2^2}{\sigma^2}\right)\,d\dvv,\end{align}
where $\dvv^\star=(0.3,0.5)$ and $\sigma=0.16$.
The advection diffusion model, defined here, can be inexpensively evaluated and therefore allows for exhaustive exploration of the behavior of our proposed multi-index collocation approach.
}


\subsection{Aero-thermo-structural model of a jet engine nozzle}
\label{sec-2-2}
\begin{figure}[htb]
\begin{center}
\includegraphics[width=\textwidth]{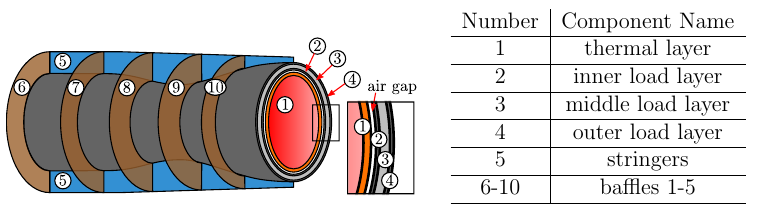}
\end{center}
\caption{(Left) schematic of the nozzle components. (Right) Component numbering. \rev{These figures can be found in \cite{Fenrich_MAA_ECCM_2018}}}
\label{fig:nozzle-components}
\end{figure}

One of the main objectives of this paper is to quantify uncertainty in a supersonic jet engine nozzle. Realistic modeling of the nozzle requires a multiphysics approach that couples aerodynamic, structural and thermal analyses. In this paper we use the model introduced in \cite{Fenrich_MAA_ECCM_2018}, and Figure \ref{fig:nozzle-components} provides a conceptual diagram of its components. This model performs a steady-state analysis of top-of-climb flight, where temperatures and stresses are the highest.
Forty random variables were identified and characterized from experimental data, simulations, or expert judgment. Each variable and its distribution is described in Appendix \ref{sec:appendix-random-variables}.

We have developed an automated suite of coupled multidisciplinary analysis tools for the steady aero-thermal-structural analysis of a supersonic nozzles called MULTI-F. MULTI-F is written in Python and calls the open source codes SU2 \cite{Economon_PCLA_AIAA_2016} and AERO-S \cite{Farhat_AEROS_2018} for fluid and thermal/structural analyses, respectively.
The internal and external flow is modeled by Euler equations which are solved using SU2. One-way feed-forward coupling is used to model the interaction between the fluid and structural analyses. Specifically the temperature obtained from the aerodynamic analysis furnishes the  boundary  condition  on  the  inner  wall  of  the  thermal  model. The  pressure  and temperatures obtained from the aerodynamic and thermal analyses then are used to provide the inner wall pressure and the temperature distribution used in a linear structural analysis. An  elastostatic  boundary value problem representing the equilibrium of internal and external forces is considered for the structural analysis while a Poisson boundary value problem representing the steady state heat transfer is considered for the thermal analysis.  The other boundary conditions are state-independent; a convection boundary condition is specified on the outer wall of the thermal model, while a fixed displacement boundary condition is imposed on the outer edges of the baffles and stringers in the structural model.

In the following we will leverage different resolutions of both the CFD and structural meshes to reduce the cost of quantifying uncertainty in four model QoI: mass, thrust, load-layer temperature failure ratio, and thermal-layer failure criteria.

\section{Multifidelity Modeling}
\label{sec-3}
\label{sec:discretization}
In this paper, we present a method to efficiently estimate statistics of the QoI~\eqref{eq:qoi-function} obtained from models of complex physical systems~\eqref{eq:govern} . For most, if not all, practical applications, the solutions of the governing equations and statistics of the resulting QoI cannot be computed analytically and numerical approximations must be employed. In this section, we discuss the impact of numerical simulation errors and response surface errors on the quantification of uncertainty.

Given a set of governing equations~\eqref{eq:govern}, we assume access to a PDE solver that approximates the solution $u$ of these equations for a given \textit{fixed} $\rvv$. We also assume this solver has a set of hyper-parameters --- mesh size, time step, maximum number of iterations, convergence tolerance, etc. --- which can be used to estimate the QoI $f$ 
at varying accuracy and cost. Differing settings for these hyper-parameters produces simulations of \textit{varying fidelities (resolution)}. We refer to approaches that leverage only one model or solver setting as {\it single-fidelity} methods and approaches that leverage multiple models and settings as {\it multi-fidelity} methods.

As stated in the introduction, quantifying uncertainty in numerical simulations can be costly, consequently it is common to construct a surrogate (function approximation) of the model parameters-to-output map. In this section, we introduce the solver parameters used to solve the two models considered in this paper (Section \ref{sec-3-1}) and then present the method used to build surrogates of a single fidelity model (Section \ref{sec-3-2}). These foundations are used to construct a multi-fidelity surrogate in Section \ref{sec-4}.


\subsection{Physical approximation}
\label{sec-3-1}
The accuracy and cost of a model is often determined by a set of solver parameters, such as mesh size, time step, tolerances of numerical solvers. Let $\na$ denote the number of such parameters for a given model. Furthermore let $l_{\ai_i} \in \mathbb{N}$ denote the number of values that each solver parameter can assume, such that the multi-index $\aindex=(\ai_1,\ldots,\ai_\na)\in\mathbb{N}^\na$, $\ai_i\in \mathcal{M}_{\ai_i}=\{1,\ldots,l_{\ai_i}\}$ can be used to index specific values of the solver parameters.
\rev{Finally let $\ffem(\rvv)$ to denote a single-fidelity physical approximation of the QoI $f(\rvv)$ using the solver parameters indexed by $\aindex$.
In this section we detail the solver indices $\aindex$ for the two models considered in this paper and describe their effect on the accuracy and cost of $\ffem(\rvv)$.}

\subsubsection{Advection-diffusion model}
\label{sec-3-1-1}
\label{sec:adv-diff-phys-discretization}
\rev{The advection-diffusion equation \eqref{eq:advection-diffusion} is solved using the Galerkin finite element method with linear finite elements and the implicit backward-Euler time-stepping scheme. Let \(h_1\) and $h_2$ denote the mesh-size in the spatial directions $\dv_1$ and $\dv_2$. Then given some base coarse discretization $h_{j,0}$, we create a sequence of $l_{\ai_j}$ uniform triangular meshes of increasing fidelity by setting \(h_{\ai_j}=h_{j,0}\cdot 2^{-\ai_j}\), $0\le\ai_j<l_{\ai_j}$, where $h_{j,0}=\frac{1}{4}$ and $l_{\ai_j}=6$, $j=1,2$. The number of vertices in the spatial mesh is $(h^{-1}_{\ai_1}+1)(h^{-1}_{\ai_2}+1)$.

The backward-Euler scheme used to evolve the advection-diffusion equation in time is an implicit method and unconditionally stable; consequently, the time-step size $\Delta t$ can be made very large. In this paper, we use an ensemble of $l_{\ai_3}$ time-step sizes \(\Delta t_{\ai_1,\ai_3}=\Delta t_{\ai_3,0}\cdot 2^{-\ai_3}\), $0\le\ai_3< l_{\ai_3}$ to solve ~\eqref{eq:advection-diffusion}  and set $l_{\ai_3}=6$ and $\Delta t_{\ai_3,0}=1/4$.
For high enough $\ai_1$ and $\ai_2$, the ensemble represents a sequence of models of increasing accuracy.

The cost of solving the advection-diffusion is dependent on both the number of mesh vertices and the number of time-steps. The number of mesh vertices is $2^{(\ai_1+2)(\ai_2+2)}$, whereas the number of time-steps is $2^{\ai_3+2}$. In the following we use a conservative estimate for the computational cost of evaluating the advection diffusion model. Specifically we assume that the cost is $O(2^{(\ai_1+2)(\ai_2+2)(\ai_3+2)})$. This cost model is representative of the best scaling that one could hope for when solving the advection diffusion on a large scale. It mimics the use of a linear solver whose costs grows linearly with the number of degrees of freedom and a finite element assembly which is constant regardless of the mesh discretization.


Figure~\ref{fig:work-vs-error} depicts the changes induced in the error of the mean of the QoI, i.e. $\mathbb{E}[f]$, as the mesh and temporal discretizations are changed.\footnote{The expectation is computed using the same 10 samples for all model resolutions. This allows the error in the statistical estimate induced by small numbers of samples to be ignored. The reference solution is obtained using the model indexed by $(6,6,6)$.} The error decreases quadratically with both $h_1$ and $h_2$ and linearly with $\Delta t$ until a saturation point is reached. These saturation points occur when the error induced by a coarse resolution in one mesh parameter dominates the others. For example the left plot shows that when refining $h_1$ the final error in $\mathbb{E}[f]$ is dictated by the error induced by using the mesh size $h_2$, provided $\Delta t$ is small enough. Similarly the right plot shows, that for fixed $h_2$, the highest accuracy that can be obtained by refining $\Delta t$ is dependent on the resolution of $h_1$. In Section \ref{sec:results-advection-diffusion}, we leverage this varying cost-vs-accuracy profile to efficiently estimate uncertainty in the advection-diffusion QoI.

\begin{figure}[htb]
\begin{center}
\includegraphics[width=\textwidth]{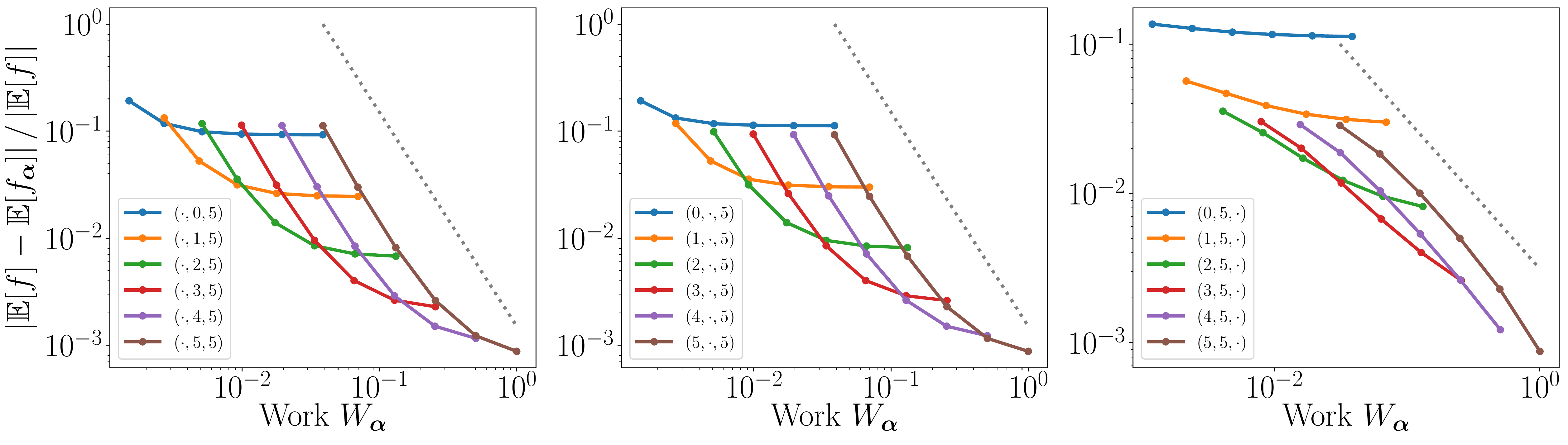}
\end{center}
\caption{\rev{The error in $\mathbb{E}[f]$ for varying mesh and temporal resolutions. The legend labels denote the mesh discretization parameter values $(\alpha_1,\alpha_2,\alpha_3)$ used to solve the advection diffusion equation. Numeric values represent discretization parameters that are held fixed while the symbol $\cdot$ denotes that the corresponding parameter is varying. The dashed lines represent the theoretical rates of the convergence of the deterministic error when refining $h_1$ (left), $h_2$ (middle), and $\Delta t$ (right).}
}
\label{fig:work-vs-error}
\end{figure}
}

\subsubsection{Aero-thermal-structural model for a jet engine nozzle}
\label{sec-3-1-2}

\begin{figure}[htb]
\begin{center}
\includegraphics[width=.5\textwidth]{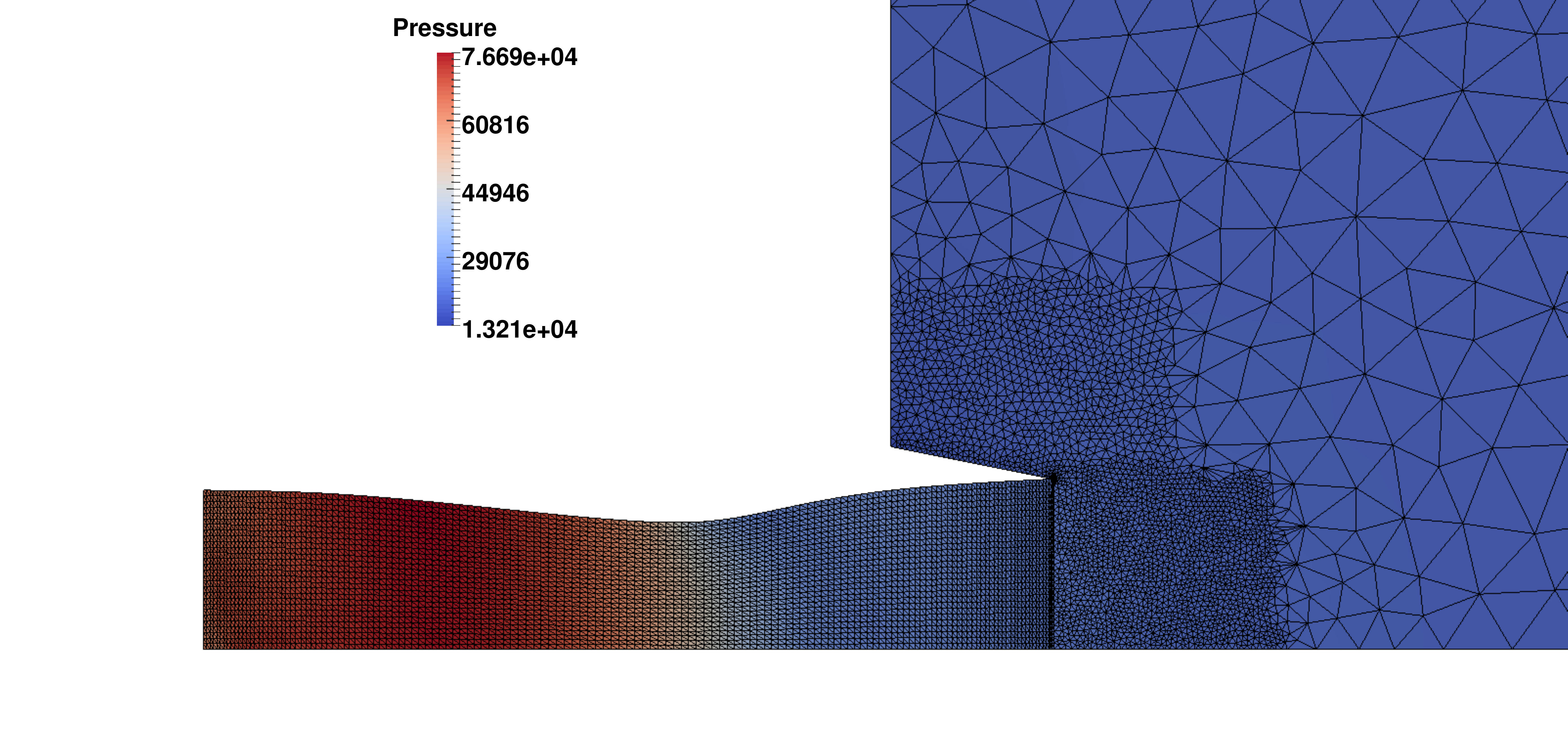}
\includegraphics[width=.49\textwidth]{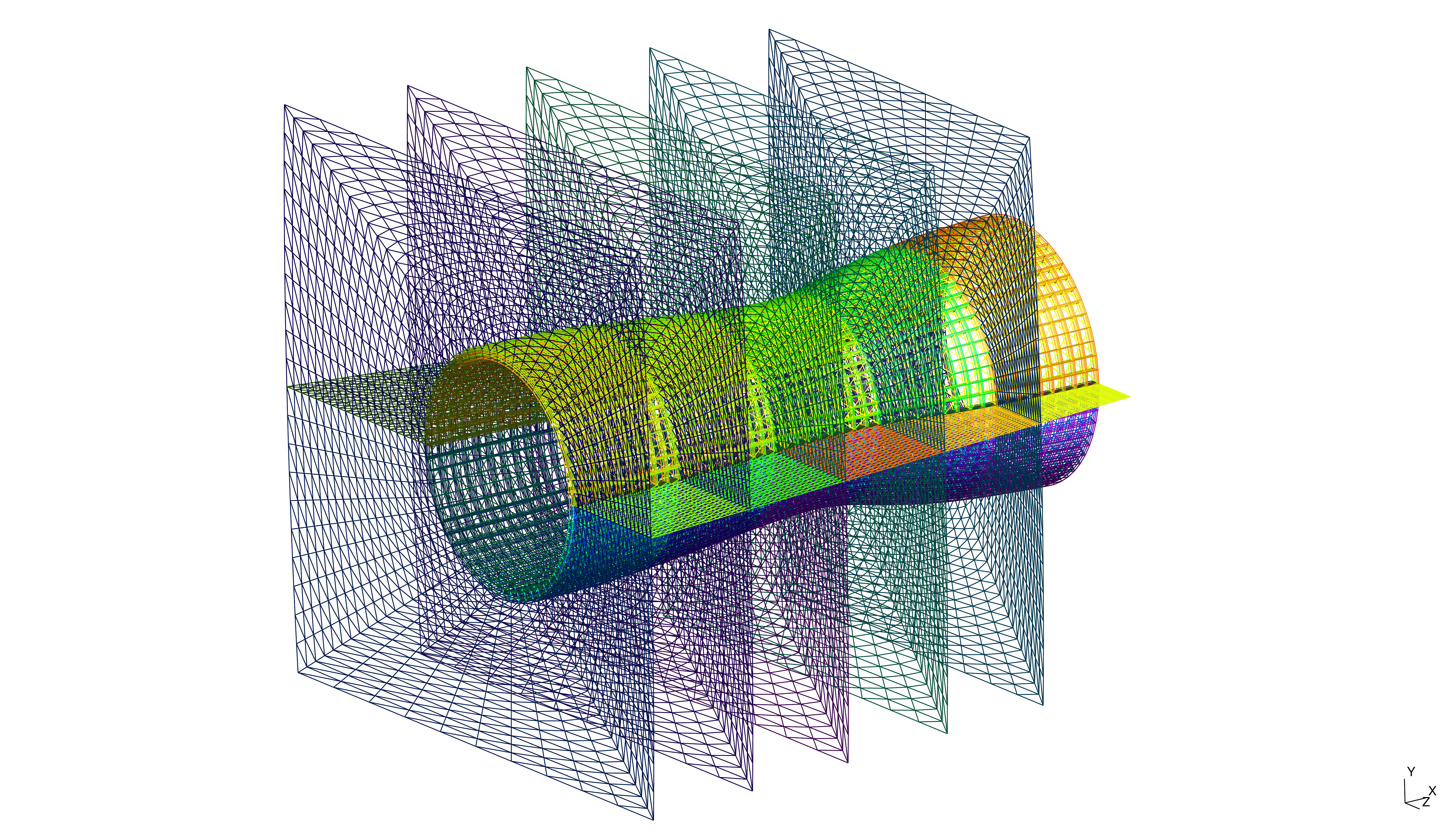}
\end{center}
\caption{(Left) The coarse ($\ai_1=0$) 2D CFD mesh for the nozzle geometry and pressure contours for the fluid moving from left to right. (Right) The coarse mesh ($\ai_2=0)$ of the structural model with geometry specified using the three-dimensional parameterization. The fluid moves from left to right.}
\label{fig:nozzle-meshes}
\end{figure}

The nozzle model has two solver parameters which control accuracy and cost. These parameters specify the mesh resolutions for the fluid and structural analyses.
The internal and external nozzle flow is calculated using SU2, an open-source software suite for multiphysics simulations \cite{Economon_PCLA_AIAA_2016}. The governing Euler equations are discretized in SU2 using a finite volume method with a standard edge-based data structure. The first discretization index $0 \le \ai_1<3$ indexes the resolution of the CFD mesh. Specifically the meshes indexed by $\aindex $ represent a  sequence  of  three  increasingly  refined meshes. The coarsest CFD mesh is depicted in the left of Figure \ref{fig:nozzle-meshes}.

The nozzle thermal and structural analyses are calculated using AERO-S, an open-source finite element method (FEM) analysis software \cite{Farhat_AEROS_2018}. The second discretization index $0\le \ai_2< 5$ of the nozzle model indexes a sequence of increasingly refined structural meshes. The coarsest structural mesh is depicted on the right of Figure \ref{fig:nozzle-meshes}. 

The relative computational cost of each nozzle model approximation is given in Table \ref{tab:nozzle-cost}. The cost is given in wall time using 4 CPUs for the CFD model and 1 CPU for the structural model (the structural model cannot currently be run in parallel).\footnote{\rev{The timings in Table~\ref{tab:nozzle-cost} were obtained with 4 Intel(R) Xeon(R) CPU E5-2640 v4 @ 2.40GHz (3.40GHz with turbo boost) processors and 64GB of RAM. The absolute numbers are not important, only the relative timings.}}
\begin{table}
\centering
\caption{Computational cost (in seconds on 4 CPUs) of the nozzle model.}
\label{tab:nozzle-cost}
\begin{tabular}[htb]{|c|c|c|c|c|c H|}
\hline 
\backslashbox{$\alpha_1$}{$\alpha_2$} &0 &1 &2 &3 &4 &6 \\ 
\hline 0 & 36.3 & 39.3 & 44.2 & 61.5 & 131.2 & 347.9 \\ 
\hline 1 & 88.5 & 90.5 & 95.7 & 113.1 & 184.2 & 386.8 \\ 
\hline 2 & 293.6 & 297.4 & 300.0 & 317.5 & 384.5 & 609.2 \\ 
\hline
\end{tabular}
\end{table}

\subsection{Stochastic approximation}
\label{sec-3-2}
Given evaluations of $\ffem(\rvv)$ at various realizations of the random parameters $\rvv$,
this section describes a method to construct surrogates that can interpolate these samples. 
%
%
%
With this goal let $\bindex=(\bi_1,\ldots,\bi_\nb)\in\mathbb{N}^\nb$ denote a multi-index that controls the computational cost of constructing a surrogate and its corresponding accuracy. For example, $\bi_i$ could be the number of univariate samples in the $i$th variable dimension, $1\le i\le \nb=d$, used to construct a tensor product interpolant. We refer to the entries of $\bindex$ as stochastic approximation hyper-parameters and denote the surrogate, built with this choice of parameters, as $\surr_{\aindex,\bindex}$. 

Numerous methods, such as those discussed in the introduction, have been developed to build surrogates of simulation models.
In this paper we choose to focus on sparse grid surrogates because the single-fidelity formulation, presented in this section, can easily be extended to use multiple model fidelities as shown in Section \ref{sec-4}.
\subsubsection{Tensor-product interpolation}
\label{sec-3-2-1}
In this section we describe tensor-product interpolation. Knowledge of tensor-product interpolation is an essential component of sparse grid approximation. Let \(\surr_{\aindex,\bindex}(\rvv)\) be an \(M\)-point tensor-product interpolant of the function \(\ffem\), with $\nb=\nv$. This interpolant is a weighted linear combination of tensor-product of univariate Lagrange polynomials 
\begin{align}
\basis_{i,j}(\rv_i) = \prod_{k=1,k\neq j}^{m_{\bi_i}}\frac{\rv_i-\unisamp{i}{k}}{\unisamp{i}{j}-\unisamp{i}{k}}, \quad i\in\{1,\ldots,d\},
\end{align}
defined on a set of univariate points $\unisamp{i}{j},j\in\{1,\ldots,m_{\bi_i}\}$.\footnote{For tensor-product approximation and sparse grid approximation we can use $\nb$ and $\nv$ interchangeably because $\nb=\nv$.} Specifically the multivariate interpolant is given by
\begin{align}
\surr_{\aindex,\bindex}(\rvv) = \sum_{\V{j}\le\bindex}
\ffem(\rvv^{(\V{j})})\prod_{i=1}^d\basis_{i,j_i}(\rv_i).
\end{align}
The partial ordering $\V{j}\le\bindex$ is true if all the component wise conditions are true.

Constructing the interpolant requires evaluating the function $\ffem$ on the grid of points
\begin{align}
\grid_\bindex = \bigotimes_{i=1}^\nv \grid_{\bi_i}^i=\begin{bmatrix}
\samp{1} & \cdots&\samp{M_\bindex}
\end{bmatrix}\in\reals^{d\times M_\bindex}
\end{align}
We denote the resulting function evaluations by 
$\values_{\aindex,\bindex}=\ffem(\grid_\bindex)=\begin{bmatrix}\ffem(\samp{1}) & \cdots& \ffem(\samp{M_\bindex})\end{bmatrix}^T\in\reals^{M_\bindex\times\nq}$, where the number of points in the grid is $M_\bindex=\prod_{i=1}^d m_{\bi_i}$.

It is often reasonable to assume that, for any $\rvv$, the cost of each simulation is constant for a given $\aindex$. So letting $W_\aindex$ denote the cost of a single simulation, we can write the total cost of evaluating the interpolant $W_{\aindex,\bindex}=W_\aindex M_\bindex$. Here we have assumed that the computational effort to compute the interpolant once data has been obtained is negligible, which is true for sufficiently expensive models $\ffem$.
In this paper, we use the nested Clenshaw-Curtis points
\begin{align}
\unisamp{i}{j}=\cos\left(\frac{(j-1)\pi}{m_{\bi_i}}\right),& & j\in\{1,\ldots,m_{\bi_i}\}
\end{align}
to define the univariate Lagrange polynomials. The number of points $m(l)$ of this rule grows exponentially with the level $l$, specifically
$m(0)=1$ and $m(l)=2^{l}+1$ for $l\geq1$.
We remark here that any univariate sequence of points can be used to construct a sparse grid. However the most efficient rules are nested and tailored to the distribution of the random variables $\rVv$. Clenshaw-Curtis points are popular for uniform variables and weighted Leja sequences \cite{Narayan_J_SISC_2014} can be used for any continuous bounded variable as well as for a broad class of unbounded variables.

When using surrogates for uncertainty quantification, it is often useful to compute the expectation of the approximation. The expectation $\mu_\bl$ of a tensor product interpolant can be computed without explicitly forming the interpolant and is given by
\begin{align}
\mu_\bl=
\myint{\sum_{\V{j}\le\bindex}
\ffem(\rvv^{(\V{j})})\prod_{i=1}^\nv\basis_{i,j_i}(\rv_i)}=\sum_{\V{j}\le\bindex} \ffem(\rvv^{(\V{j})}) v_{\V{j}}.
\end{align}
The expectation is simply the weighted sum of the Cartesian-product of the univariate quadrature weights $v_{\V{j}}=\prod_{i=1}^\nv\int_{\rvvsupp_i}{\basis_{i,j_i}(\rv_i)}\dx{\pbwt(\rv_i)}$, which can be computed analytically.

\subsubsection{Sparse grid interpolation}
\label{sec-3-2-2}
\label{sec:sparse-grid}
The computational effort required to generate function values for a tensor-product interpolation is proportional to the number of grid points, which grows exponentially with the number of random parameters $\nv$. In this section we introduce sparse grid interpolation which can construct multivariate function approximations with accuracy comparable to that of tensor-product approximations, but at a fraction of the cost.

The sparse grid interpolant and expectation of the function $\ffem$ can be expressed as a linear combination of low-resolution tensor-product interpolants \cite{Barthelmann_NR_ACM_2000}
\begin{align}
\label{eq:general-smolyak}
\sgrid(\rvv)
=\sum_{\bindex\in\iset}c_{\bindex}\,\surr_{\aindex,\bindex}(\rvv) & & 
\mu_{\aindex,\iset}=\sum_{\bindex\in\iset}c_{\bindex}\,\mu_{\aindex,\bindex},
\end{align}
\rev{where the coefficients $c_{\bindex}$ are given by
\begin{align}
c_{\bindex} = \sum_{\V{j}\in\otimes_{i=1}^\nv\{0,1\}}\left(-1\right)^{\lvert \V{j}\rvert_1}\bigchi_\iset(\bindex+\V{j}), && \bigchi_\iset(\V{k})=
\begin{cases}
  1 & \V{k} \in \iset\\
  0 & \text{otherwise}
\end{cases}. 
\end{align}}

Recalling that $\grid_\bindex$ are the grid points needed to build the tensor product interpolation $\surr_{\aindex,\bindex}$, then the points in a sparse grid are
$
\grid_{\iset} = \bigcup_{\bindex\in\iset} \grid_\bindex.
$
When using nested univariate sequences to construct the tensor product grids we can write
\begin{align*}
  \grid_{\iset}=\bigcup_{\bindex\in\iset}\diffgrid_\bindex, & & \diffgrid_\bindex=\grid_\bindex\setminus\bigcup_{\substack{\bindex^\star\in\iset\\ \bindex^\star\le\bindex}} \grid_{\bindex^\star},\quad \norm{\bindex}{1}>0.
\end{align*}
where $\diffgrid_\V{0}=\emptyset$ and $\diffgrid_\bindex\cap\diffgrid_{\bindex^\star}=\emptyset$ if $\bindex\neq\bindex^\star$. This decomposition of the sparse grid points indicates that many sparse grid points are shared among its constituent tensor-product grids, and thus the total number of points in the sparse grid satisfies $$M(\iset)=\mathrm{card}\left(\grid_{\iset}\right)=\mathrm{card}\left(\bigcup_{\bindex\in\iset}\grid_\bindex\right) \le \sum_{\bindex\in\iset}M_\bindex.$$

The index set $\indexset$ 
used in a sparse grid approximation \eqref{eq:general-smolyak} can be tailored to minimize error whilst also minimizing cost, i.e. the number of function evaluations. The only restriction on the index set is that it must be downward closed. An index set $\indexset$ is downward closed if $\bindex_1\in\indexset$ implies that $\bindex_2\in\indexset$  for all $\bindex_2\le\bindex_1$.

Finding an efficient index set $\iset$ can be cast as an optimization problem. With this goal, let
\begin{align}
\derr=\norm{\surr_{\aindex,\iset\bigcup\{\bindex\}} -\sgrid}{} & &
\dwork=\abs{\mathrm{Work}[\surr_{\aindex,\iset\bigcup\{\bindex\}}]-\mathrm{Work}[\sgrid]}
\end{align}
respectively denote the norm of the difference between the sparse grid with and without the interpolant $f_{\aindex,\bindex}$ and the work needed to add $f_{\aindex,\bindex}$ to the sparse grid. 
Noting that the error in the sparse grid \rev{satisfies $\norm{f-\surr_{\aindex,\iset}}{} \le \sum_{\bindex\notin\iset} \derr$}, we can formulate finding a quasi-optimal index set as a binary knapsack problem \cite{Bungartz_G_AN_2004}
\begin{align}\label{eq:knapsack}
  \max \sum_{[\bindex]\in \mathbb{N}_0^{\nb}} \derr \delta_{\bindex} \quad \text{such that} \quad \sum_{[\bindex]\in \mathbb{N}_0^{\nb}} \dwork \delta_{\bindex} \le W_{\max}, \quad\delta_{\bindex}\in\{0,1\}
\end{align}
for some maximum work $W_{\max}$. The solution $\iset=\{\bindex\mid \delta_{\bindex}=1\}$ to this problem balances the computational work of adding a specific interpolant with the reduction in error that would be achieved.\footnote{\rev{Note $\iset$ is only ``quasi'' optimal since the error decomposition maximized in the knapsack problem is only an upper bound for the true error.}}

One approach for solving the knapsack problem is to use estimates of the error $\derr$ derived from assumptions on the smoothness of the function $f_\aindex$. The resulting grid depends on the space of functions that $f_\aindex$ belongs to, and not on the function itself. Assuming that the function $f_\aindex$ has continuous mixed derivatives of order $r$ in each random variable dimension and $\dwork=\mathrm{card}(\diffgrid_\bindex)$, yields the traditional isotropic sparse grid index set and corresponding sparse grid coefficients given below

\begin{align}
  \label{eq:smolyak-coeffs}
  \iset(l) = \{\bindex \mid (\max(0,l-1)\le\lVert\bindex\rVert_1\le l-\nv-2\} & &
c_{\bindex}=\left(-1\right)^{l-\lVert\bindex\rVert_1}{\nv-1 \choose l-\lVert\bindex\rVert_1}.
\end{align}
Note again we have used the fact $\nb=\nv$ and that the form of the isotropic index set and coefficients reflect our choice to set $\bi_i\ge 0$.\footnote{Much of the literature assumes $\bi_i>0$, in which case \(\iset(l) = \{\bindex \mid (l+1\le\lVert\bindex\rVert_1\le l+\nv\}\) and $c_{\bindex}=\left(-1\right)^{l+\nv-\lVert\bindex\rVert_1}{\nv-1 \choose l+\nv-\lVert\bindex\rVert_1}.$}  The sparse grid using \eqref{eq:smolyak-coeffs} is often referred to as the Smolyak algorithm \cite{Smolyak_SMD_1963}.

There is no exact formula for the number of points in an isotropic sparse grid. However, for fixed level $l$  and large enough $\nv$, the number of points in the Clenshaw-Curtis isotropic sparse grid satisfies
$$M_{\iset(l)}\approx \frac{2^l}{l!}\nv^l,$$ where we use $\approx$ to denote the strong equivalence of sequences \cite{Novak_R_CA_1999}. This is much smaller than the number of points in the tensor-product grid
\begin{align*}
  M_{\V{l}}=m_{l}^\nv=\left(2^{l}+1\right)^d & & \V{l}=(l,\ldots,l)
\end{align*}

In Table~\ref{tab:sparse-grid-growth}, we compare the number of points $M_{\iset(l)}$ in a isotropic sparse grid and the number of points $M_\V{l}$ in a tensor-product grid for increasing number of random parameters and fixed $l=2$. It is clear that the number of points in the sparse grid increases much more slowly than the number of points in the tensor-product grid.
\begin{table}
  \centering
\begin{tabular}[htb]{|c|c|c|c|c|c|}
\hline  {$d$} &2 &5 &10 &20 &40 \\ 
\hline $M_{\mathcal{I}(l)}$ & 2.5e+01 & 3.1e+03 & 9.8e+06 & 9.5e+13 & 9.1e+27 \\ 
\hline $M_{\boldsymbol{\beta}}$ & 1.3e+01 & 6.1e+01 & 2.2e+02 & 8.4e+02 & 3.3e+03\\
\hline
\end{tabular}
\caption{The number of points $M_{\iset(l)}$ in a isotropic sparse grid and the number of points $M_\V{l}$ in a tensor-product for increasing number of random parameters and fixed $l=2$.}
\label{tab:sparse-grid-growth}
\end{table}

For a function with $r$ continuous mixed derivatives, the isotropic level-$l$ sparse-grid based on Clenshaw-Curtis abscissas with $M_{\iset(l)}$ points satisfies \rev{\cite{Barthelmann_NR_ACM_2000}}:
\begin{equation}\label{eq:sparse-grid-error}
\norm{f-\surr_{\iset(l)}}{L^\infty(\rvvsupp)} \le C_{\nv,r} M_{\iset(l)}^{-r}(\log M_{\iset(l)})^{(r+2)(\nv-1)+1}
\end{equation}
and the tensor-product interpolant satisfies
\begin{equation}\label{eq:tp-error}
\norm{f-\surr_\V{l}}{L^\infty(\rvvsupp)} \le K_{\nv,r} M_{\V{l}}^{-r/d}
\end{equation}
where the constants $C_{\nv,r}$ and  $K_{\nv,r}$ depend on the smoothness $r$ and the dimension $\nv$, but not on the number of points in the grid. From \eqref{eq:sparse-grid-error}, we can see that the error in a sparse grid approximation depends only weakly on the dimension $\nv$. In contrast, the error in a tensor product interpolant \eqref{eq:tp-error} depends exponentially on the number of random parameters $\nv$. For both sparse grid and tensor product interpolation, the error in each approximation improves as the smoothness $r$ increases.

We remark here that, although an improvement over tensor-product interpolation, the logarithmic term in the sparse grid error bound can become large when the number of model variables becomes sufficiently large. In these situations, adaptive strategies for constructing the index set $\iset$ have proven to be very effective. We present one such adaptive strategy in Section \ref{sec:adaptivity}.


\subsubsection{An example of isotropic sparse grid approximation}
\begin{figure}[ht]
\begin{center}
\includegraphics[width=\textwidth]{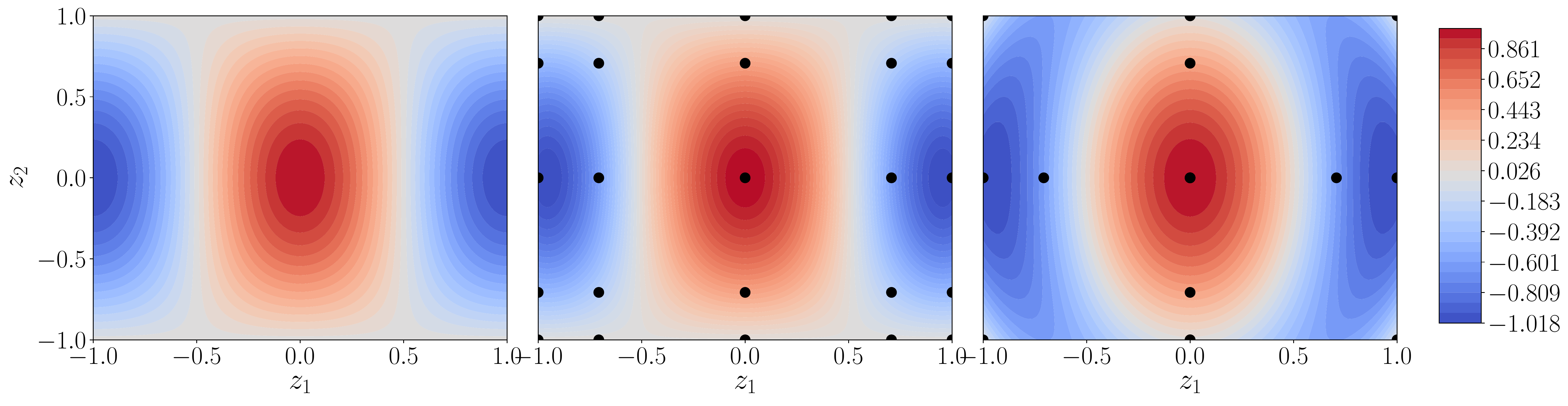}
\end{center}
\caption{(Left) the function $\ffem(\rvv)=\cos(2\pi\rv_1)\cos(\pi\rv_2)$. (Middle) the tensor product interpolant $\surr_{\aindex,(2,2)}(\rvv)$. (Right) the isotropic sparse grid interpolant $\surr_{\aindex,\iset(2)}(\rvv)$. The black circles depict the samples used to build each interpolant.}
\label{fig:tensor-product-approx}
\end{figure}

In this section we use a simple example to provide insight into the construction and utility of sparse grid approximations.
Consider the function\footnote{\rev{Note the function in \eqref{eq:cosine} has no physical solver parameters, however we keep the subscript $\aindex$ to be consistent with the notation in the rest of the paper.}}
\begin{align}\label{eq:cosine}\ffem(\rvv)=\cos(2\pi\rv_1)\cos(\pi\rv_2)\end{align}
which depends on two random variables $\rVv=(\rV_1,\rV_2)$ which are uniformly distributed in $\rvvsupp=[-1,1]^2$. The contours of this function are depicted in Figure~\ref{fig:tensor-product-approx}. Figure~\ref{fig:tensor-product-approx} also plots the points used to build a $l=2$ tensor-product approximation, the points used to construct a $l=2$ sparse grid approximation, and the contours of the resulting approximations.

\begin{figure}[ht]
\begin{center}
\includegraphics[width=\textwidth]{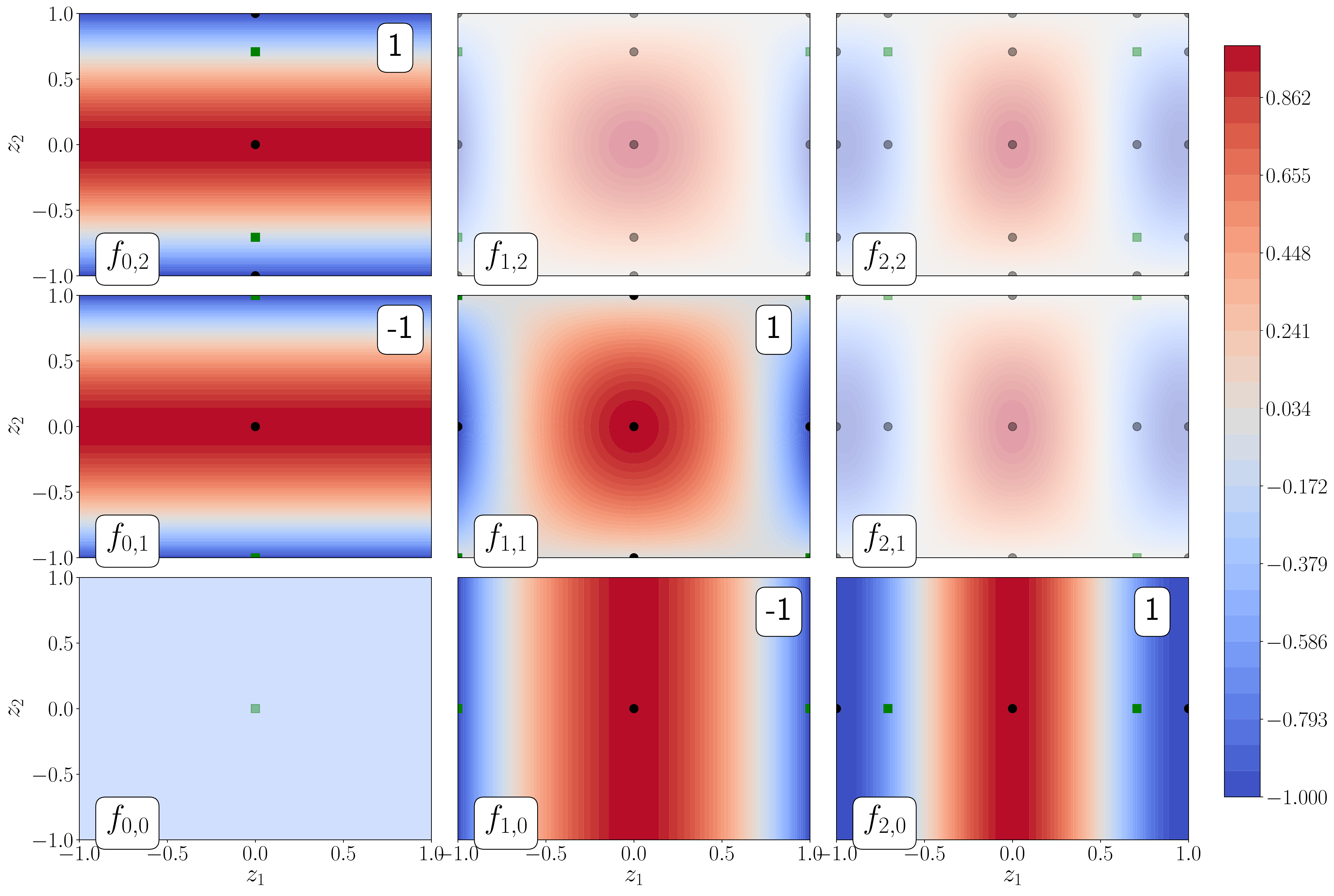}
\end{center}
\caption{The tensor-product grids and contours of the interpolants used to construct the level $l=2$ isotropic sparse grid $\isgrid{2}(\rvv)$. The numbers in the top right of each plot are the coefficients $c_\bindex$ in \eqref{eq:smolyak-coeffs}. \rev{The tensor interpolants without coefficients do not contribute to the sparse interpolant}. The green squares depict the unique samples $\diffgrid_\bindex$ of the tensor-product grid $\grid_\bindex$ that do not exist in another tensor-product grid $\grid_{\bindex^\star}$, $\bindex^\star\le\bindex$. Both the samples depicted by the black circles and the green squares are used to build an interpolant.}
\label{fig:tensor-product-spaces}
\end{figure}

The isotropic sparse grid, depicted in the right plot of Figure ~\ref{fig:tensor-product-approx}, is the weighted linear combination of the tensor-product interpolants shown in Figure \ref{fig:tensor-product-spaces}. The sparse grid points are the union of the points of the low resolution grids with non-zero coefficients $c_\bindex$. The sparse grid achieves a similar accuracy to the tensor-product interpolant $\surr_{\aindex,(2,2)}(\rvv)$ (see the middle plot of Figure ~\ref{fig:tensor-product-approx} and the top right plot of Figure \ref{fig:tensor-product-spaces}) but does not need to use the unique points $\diffgrid_\bindex$ (depicted by green squares) in the finest tensor-product grids (indexed by $(1,2),(2,1)$ and $(2,2)$). Moreover despite the sparse grid being comprised of low-resolution/low-accuracy tensor-product interpolants, the linear combination of these interpolants produces a much more accurate approximation than any single low-resolution interpolant. Avoiding the construction of fine tensor-product interpolants, which require large numbers of points, is the primary mechanism by which sparse grids \rev{mitigate} the curse of dimensionality.

Note the approximation $\surr_{0,0}$ is not included in the sparse grid interpolant. As the level $l$ is increased more lower-level tensor-product interpolants will be excluded from the sparse grid approximation.


\section{Multi-index collocation}
\label{sec-4}
\label{sec:misc}

Sparse-grids and other surrogate methods can \rev{mitigate} the curse of dimensionality by exploiting structure in the function being approximated. However, many if not all surrogate methods can be infeasible in high-dimensions if the computational cost of evaluating a model is high. In this section, we present a method for building surrogates in such situations, which leverages simulations of lower-fidelity models of reduced cost to increase the tractability of building surrogates of a high-fidelity model. We focus on sparse grid function approximation; however, in Section~\ref{sec:sensitivity-analysis}, we briefly discuss how the algorithm we use for managing evaluations of models of varying fidelity can be generalized to other types of approximations\rev{, such as polynomial chaos expansion constructed on arbitrary points sets, Gaussian processes, low-rank approximations, or even neural networks.}

Using the triangle inequality, it is straightforward to decompose the error in the sparse grid approximation $\sgrid$ of the exact model $f$, built using evaluations of $\ffem$, into the sum of the a deterministic and stochastic approximation errors.
\begin{align}
\label{eq:error-break-down}
\mynorm{f-\sgrid}\le\mynorm{f-\ffem}+\mynorm{\ffem-\sgrid}
\end{align}
Here the first term on the right-hand side represents the deterministic error and the second term represents the stochastic error. These two errors should be balanced to minimize inefficiency resulting from over-resolving a non-dominant source of error.
\rev{Figure \ref{fig:work-vs-error} demonstrates that mesh discretizations can effect the accuracy of statistics computed using a numerical model. For example, when the time step is to large, refining the spatial mesh has no effect on error in the mean $\mathbb{E}[\ffem]$. Analogously \eqref{eq:error-break-down} implies that refining the mesh when the stochastic error is larger will have no effect on overall error $\mynorm{f-\sgrid}$.}


Multi-level/Multi-index collocations methods \cite{DeBaar_R_PAG_2017,HajiAli_NTT_CMAME_2016,Ng_E_AIAA_2012,Teckentrup_JWG_SIAMUQ_2015} combine a sequence of approximate physical models with a sequence of interpolants both to balance physical and stochastic errors and to reduce the computational cost of achieving a specified level of accuracy. In this paper, we \rev{extend} the multi-index stochastic collocation method (MISC) presented in \cite{HajiAli_NTT_CMAME_2016}. The other \rev{multi-level} methods cited can be formulated as specialized instances of MISC.

Given an ensemble of models \rev{$\mathcal{M}_{\aindex}=\otimes_{i=1}^\na \mathcal{M}_{\ai_i}$} of varying physical fidelity \rev{(numerical discretization)}, which can be used to approximate the QoI $f$, MISC approximates the QoI and its mean with the following generalization of Equation~\eqref{eq:general-smolyak}
\begin{align}
\label{eq:general-misc}
\misc(\rvv)
=\sum_{[\aindex,\bindex]\in\jset}c_{\aindex,\bindex}\,\surr_{\aindex,\bindex}(\rvv), & & 
\mu_{\jset}=\sum_{[\aindex,\bindex]\in\jset}c_{\aindex,\bindex}\,\mu_{\aindex,\bindex}.
\end{align}
MISC works most effectively provided the following assumption
\begin{assumption}
Given a sequence of tolerance parameters, 
the error in the successive approximations of $f$ decreases as the fidelity increases, i.e. $\norm{\surr_{\aindex^\star}-f}{L^p(\rvvsupp)}<\norm{\surr_{\aindex}-f}{L^p(\rvvsupp)}$, if $\aindex^\star > \aindex$.
\end{assumption}

Similarly to sparse grid function approximation (see Section \ref{sec:sparse-grid}), the structure of the index set $\jset$ of the MISC approximation \eqref{eq:general-misc} directly controls the computational cost of building, and the accuracy of, a multi-index surrogate. However, unlike the sparse grid index set $\iset$, the MISC index set $\jset$ not only controls the allocation of samples to each random variable, but also the allocation of resources to each of the model discretizations.

Given a MISC approximation $\misc$, let $\surr_{\jset\bigcup\{[\aindex,\bindex]\}}$ be the new MISC approximation obtained by including a new tensor-product interpolant $\surr_{\aindex,\bindex}$. Similarly to Section~\ref{sec:sparse-grid}, let
\begin{align}
\derr=\norm{\surr_{\jset\bigcup\{[\aindex,\bindex]\}} -\misc}{} & &
\dwork=\abs{\mathrm{Work}[\surr_{\jset\bigcup\{[\aindex,\bindex]\}}]-\mathrm{Work}[\misc]}
\end{align}
denote the difference between two successive MISC approximations and the work needed to update the MISC approximation, respectively. For a given $\aindex$, the cost of running a simulation at a single realization of $\rvv$ is again denoted as $W_\aindex$, such that we have $\dwork = W_\aindex\mathrm{card}(\diffgrid_{\aindex,\bindex})$.

With the goal of generating a MISC approximation that balances the computational work and error, we can again formulate and solve a knapsack problem to define the index set $\jset$. The knapsack problem is identical to \eqref{eq:knapsack}, except now the optimization is not just over $\bindex\in\mathbb{N}$, but also over $\aindex\in\mathbb{N}$. If the quantities $\derr$ and $\dwork$ are available, the quasi-optimal\footnote{The solution is only quasi optimal because it is based upon an upper bound on the approximation error, that is $\norm{f-\misc}{}\le \sum_{[\aindex,\bindex]\notin \jset} \derr$.} solution to the knapsack problem is
\begin{align}
\jset=\left\{[\aindex,\bindex]\in\mathbb{N}_0^{\na+\nb} \mid \frac{\derr}{\dwork}>\epsilon \right\},
\end{align}
where $\epsilon$ is chosen to reflect the desired accuracy in the MISC approximation \cite{HajiAli_NTT_CMAME_2016}. Computing this index set is not practical as it requires computing all possible interpolants and then choosing the approximations which significantly contribute to the MISC error. The interpolants ignored would contribute little to error but may incur significant cost.  In this paper, we propose using an adaptive algorithm to greedily choose the most cost effective interpolants to add to a index set $\jset$. We will present our adaptive algorithm in Section \ref{sec:adaptivity}; however, we first provide an example of using MISC with an isotropic index set to approximate a simple algebraic ensemble of models. \rev{This example is intended to develop intuition about the fundamental mechanisms which lead to the improved efficiency of MISC, as compared to single fidelity strategies.}

\subsection{An example of isotropic MISC approximation}
\begin{figure}[ht]
\begin{center}
\includegraphics[align=c,width=\textwidth]{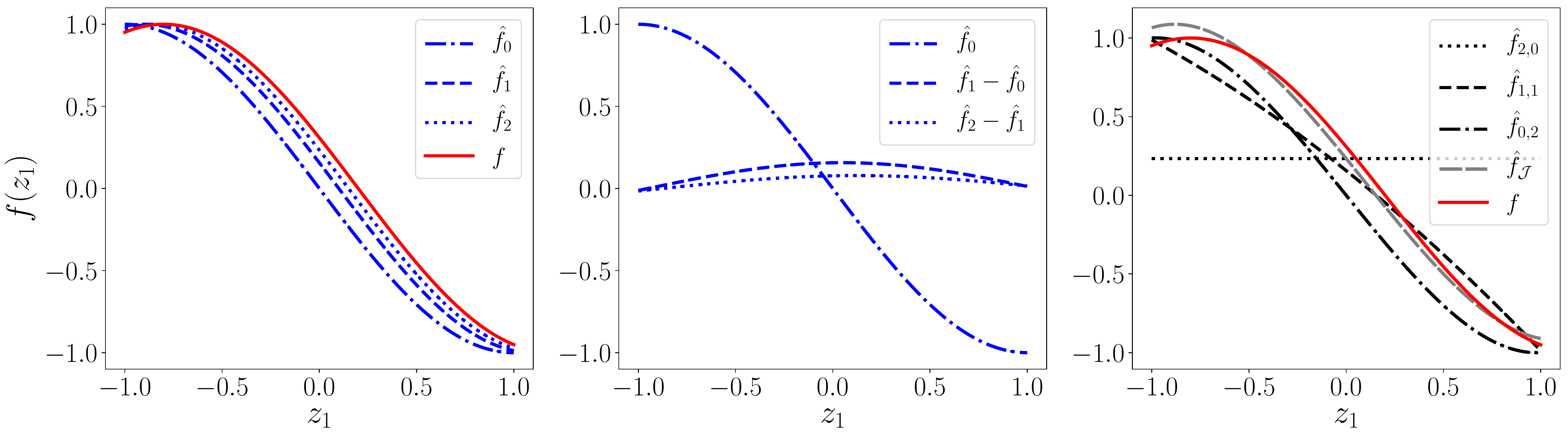}
\end{center}
\caption{(Left) The function $f$ \eqref{eq:misc-example} and the approximations $\surr_{\ai_1}$. (Middle) The discrepancies between successive approximations $\surr_{\ai_1}$ and $\surr_{\ai_1+1}$. (Right)
  The most accurate (largest $\bi_1$) single fidelity interpolants $f_{\ai_1,\bi_1}$, \rev{used by MISC}, which use only evaluations of one approximate model $\surr_{\ai_1}, \ai_1=0,1,2$, and the MISC approximation $\misc$ which combines each of these interpolants (and some less accurate interpolants).}
\label{fig:misc-discrepancies}
\end{figure}

The isotropic MISC approximation of a function $f$ is
\begin{align}
\label{eq:isotropic-misc}
\imisc{l}=\sum_{[\aindex,\bindex]\in\jset(l)} \left(-1\right)^{l-\lVert\aindex+\bindex\rVert_1}{\na+\nb-1 \choose l-\lVert\aindex+\bindex\rVert_1}\surr_{\aindex,\bindex}
\end{align}
where \(\jset(l) = \{[\aindex,\bindex] \mid (\max(0,l-1)\le\lVert\aindex+\bindex\rVert_1\le l+\na+\nb-2\}\). Consider the MISC approximation of the function
\begin{equation}\label{eq:misc-example}
\rev{f(\rv)=\cos(\frac{\pi}{2}(\rv_1+\frac{4}{5}))}
\end{equation}
using the sequence of approximations \rev{$\surr_{\ai_1}=\cos(\frac{\pi}{2}(\rv_1+\frac{4}{5}+\epsilon_{\ai_1}))$} where $\epsilon_{i}>\epsilon_{{i+1}}\ge 0$. These approximations, shown in the left of Figure \ref{fig:misc-discrepancies}
for \rev{$\epsilon_{\ai_1=0}=1/5$, $\epsilon_{\ai_1=1}=1/10$ and $\epsilon_{\ai_1=2}=1/20$,}
converge to the true function as epsilon is driven to zero, i.e. \(\surr_{\ai_1}(\rvv)\rightarrow f(\rvv)\) as \({\ai_1}\rightarrow \infty\).

\begin{figure}[ht]
\begin{center}
\includegraphics[align=c,width=0.95\textwidth]{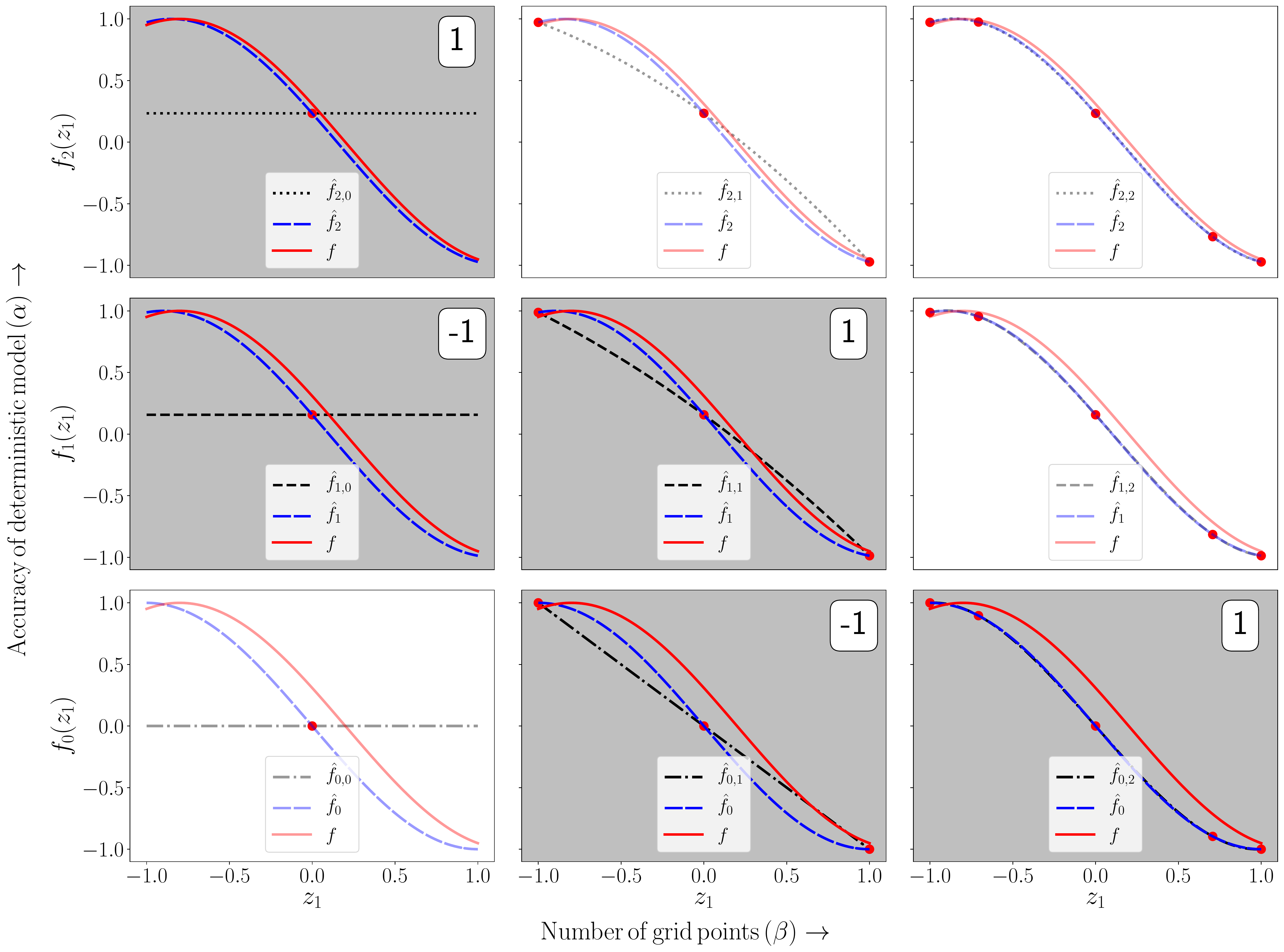}
\end{center}
\caption{The tensor-product grids and interpolants of the level $l=2$ isotropic MISC approximation $\surr_{\jset(2)}(\rvv)$.  The grids plotted with a gray background are included in the MISC approximation, whereas the remaining grids are ignored. 
The numbers in the top-right of each plot are the coefficients $c_{\aindex,\bindex}=\left(-1\right)^{l-\lVert\aindex+\bindex\rVert_1}{\nv-1 \choose l-\lVert\aindex+\bindex\rVert_1}$ in \eqref{eq:isotropic-misc}. The dash-dotted black lines represent the tensor-product interpolants $\surr_{\aindex,\bindex}$. These interpolants converge to the lower fidelity models $\surr_{\aindex}$ (dashed blue lines) as $\bindex$ increases, and to the high-fidelity model $f$ (solid red line) as both $\aindex$ and $\bindex$ are increased.}
\label{fig:misc-sparse-grid-tensor-products}
\end{figure}

In Figure \ref{fig:misc-sparse-grid-tensor-products} we plot the tensor-product grids and corresponding interpolants that make up the MISC approximation. The interpolants $\surr_{{\ai_1},\bindex}$ converge to $\surr_{{\ai_1}}$ as $\bindex$ increases. However for low values of ${\ai_1}$ the deterministic prediction error is large. Refinement of both ${\ai_1}$ and $\bindex$ is needed for the interpolants to converge to $f$.

Let $\kset=\{[\aindex^\star,\bindex^\star] \mid \bindex^\star\in\jset,\aindex^\star=\aindex\}$ denote the set of unique stochastic multi-indices for a given deterministic refinement $\aindex$. Then for each physical discretization $\aindex$ the MISC approximation requires evaluating the function $\ffem$ on the points
\begin{align}
\grid_{\kset}=\bigcup_{\bindex\in\kset}\grid_{\bindex}
\end{align}
where $M_\kset=\mathrm{card}\left(\grid_{\kset}\right)$ is the number of samples, of the random variables, at which $\ffem$ is evaluated.
The total number of samples used by the MISC approximation, for any index set $\jset$, is
$$
M_\jset=\sum_{\aindex \in \mathcal{L}}M_\kset
$$
where $\mathcal{L}=\{\aindex\mid\aindex\in\jset\}$ is the set of unique physical discretizations indices.

As can be seen in Figure \ref{fig:misc-sparse-grid-tensor-products}, the number of evaluations allocated to an approximate model $f_{\ai_1}$ decreases as the level of fidelity $\ai_1$ increases, i.e. $M_\kset>M_{\mathcal{K}_{\jset\mid\aindex^\star}}$, if $\aindex^\star>\aindex$.
For example, 5 samples are allocated to evaluating $\surr_0$, 3 to $\surr_1$ and 1 to $\surr_2$.

Despite the fact only one sample is allocated to evaluating $\surr_2$, the MISC approximation, shown in the right plot of Figure \ref{fig:misc-discrepancies}, is more accurate than any of the tensor-product interpolants that make up the MISC approximation (grids with a gray background). Understanding of this observation can be gained by viewing the discrepancies between models $\surr_{{i+1}}-\surr_{i}$ depicted in the middle plot of Figure  \ref{fig:misc-discrepancies}. The magnitude of these discrepancies decreases substantially as $\ai_1$ increases. Thus a constant approximation $\surr_{2,0}$ of $\surr_2$ is sufficient to balance the deterministic and stochastic errors, \rev{i.e. the error caused by non-zero values of $\epsilon$ and using a finite number of samples to build the interpolants, respectively}. This constant approximation is useful because it reduces the bias in the MISC interpolant by decreasing the physical discretization error. In contrast the approximation $\surr_{0,2}$ is useful as it captures the variation of the function $f$ with respect to the variable $\rv_1$ thereby reducing the stochastic error.

We remark here that the total number of samples used to build the MISC approximation is much larger than the number of samples used to build $\surr_{2,2}$. In general, the number of points in a MISC approximation is much greater than the number of points in a single fidelity sparse grid approximation of the highest fidelity model $f_{\aindex^\star}$, with equivalent cost, i.e.  $M(\jset)\gg M_{\aindex^\star,\iset}$. However, for non-trivial models, the total work needed to build the MISC approximation $\surr_\jset$ will be less than the work needed to build the single fidelity sparse grid approximation $f_{\aindex^\star,\iset}$ of $f_{\aindex^\star}$, since many of these samples occur at lower fidelity.  

\subsection{Adaptivity}
\label{sec-4-1}
\label{sec:adaptivity}
The isotropic MISC approximation \eqref{eq:isotropic-misc} balances computational expense with accuracy. The same number of model evaluations are used for two physical discretization parameters $\aindex\neq\aindex^\star$ if $\norm{\aindex}{1}=\norm{\aindex^\star}{1}$. For example the same number of model evaluations will be used to evaluate $\ffem$ and $\surr_{\aindex^\star}$ if $\na=2$ and $\aindex=(2,1)$ and $\aindex^\star=(1,2)$. The isotropic formulation does not account for differences in the cost of evaluating the two models. If $\ffem$ is twice as expensive as evaluating $\surr_{\aindex^\star}$, then ideally this should be reflected in the sample allocation used by the MISC approximation. The isotropic formulation also assumes each variable and $\surr_{\aindex^\star}$ contributes equally to the error in the MISC approximation. However, often only a subset of variables significantly influences the variability of a function. This parameter sensitivity should also be accounted for when building the MISC approximation.


We propose using a greedy algorithm to adaptively choose the most cost effective interpolants to add to the MISC approximation. Our adaptive algorithm is based on the algorithm proposed in \cite{Gerstner_G_C_2003} and pseudo-code for the algorithm 
is shown in Algorithm \ref{alg:dim-adaptivity}. The algorithm is initialized with a MISC approximation comprising of a set of indices $\jset$ representing the current approximation, a set of active indices $\cA$ that indicate the physical discretization, random dimensions for potential refinement, and sparse grid coefficients $\mathcal{C}=\{c_\bl\}_{\bl\in\jset}$. Often $\cL=\{\V{0}\}$, $\cA=\{{\V{e}_k},k\in[\na+\nb]\}$, and $\mathcal{C}=\{c_{0,0}=1\}$. Here we use $\V{e}_k$ to denote the unit vector with a value of 1 in the $k$th position.

\begin{algorithm}
\caption{\texttt{INTERPOLATE}[$\{\ffem(\rvv)\}$,$\mathcal{C}$,$\cL$,$\cA$,$\tau$,$W_{\max}$]$\rightarrow \misc$}
\label{alg:dim-adaptivity}
\begin{algorithmic}[1]
\While {NOT \texttt{TERMINATE}[$\cA$,$N$,$\tau$,$W_{\max}$]}
 \State $\blsb:= \argmax_{\blb\in\cA} \gamma_\bl$ \Comment{Determine the index with the highest priority}
 \State $\cA:=\cA\setminus \blsb$ \Comment{ Remove $\blsb$ from the active set}
 \State $\mathcal{C}$=\texttt{UPDATE}[$\jset,\blsb,\mathcal{C}$]   \Comment{Update the sparse grid coefficients}
 \State $\cL:= \cL \cup \blsb$   \Comment{ Add $\blsb$ to the MISC set}
 \State $\mathcal{R}:=$\texttt{REFINE}[$\bls$,$\cL$] \Comment{ Find all admissible forward neighbors of $\blsb$}
 \State $\gamma_{\bl}:= $ \texttt{INDICATOR}[$\bl$]$\;\forall\; \blb\in\mathcal{R}$    \Comment{ Calculate the priority of the neighbors}
 \State $\cA:= \cA \cup \mathcal{R}$ \Comment{ Add the forward neighbors to the active index set}
\EndWhile
\end{algorithmic}
\end{algorithm}

The algorithm begins by selecting the index $\blsb \in \cA$ with the largest refinement indicator $\gamma_\bl$. The interpolant $\surr_{\bls}$ is added to the sparse grid and 
its indices are used to identify new candidate interpolants for consideration in reducing
the approximation error. This process continues until a computational budget ($W_{\max}$) limiting the total computational work ($\sum_\jset \dwork$) is reached or until a global error indicator (e.g., $\sum_\cA \derr$) drops below a predefined threshold $\tau$. These exit criteria are checked using the \texttt{TERMINATE} routine in step 1 of Algorithm \ref{alg:dim-adaptivity}. The \texttt{UPDATE} function from step 4 of Algorithm \ref{alg:dim-adaptivity} is summarized in Algorithm \ref{alg:coeff-update}; it adds the interpolant $\surr_{\bls}$ to the MISC approximation. The \texttt{INDICATOR} routine from step 7 of Algorithm \ref{alg:dim-adaptivity} controls which indices are added to the sparse grid via the use of index error and global error metrics. These indicators provide estimates of the contribution of an index to reducing the error in the interpolant, and of the error in the entire interpolant, respectively.

\begin{algorithm}
\caption{\texttt{UPDATE}[$\mathcal{C}_\jset,\cL,\blsb]\rightarrow \mathcal{C}_{\jset\bigcup\{\blsb\}}$}
\label{alg:coeff-update}
\begin{algorithmic}[1]
\State $c_\bls=1$
 \For {$\blb \in \jset$}
  \State $\Delta = \blb-\blsb$
     \If {$\min_{i\in[\na+\nb]} \Delta \ge 0$  \textbf{and} $\max_{i\in[\na+\nb]} \Delta \le 1$}
         \State $c_\bl=c_\bl+(-1)^{\norm{\Delta}{1}}$
     \EndIf
\EndFor
\end{algorithmic}
\end{algorithm}

In this paper we use the following indicators to guide refinement
\begin{align}
\label{eq:dim-surplus-indicator}
  \derr^\mu = \frac{1}{\left\lvert \surr_{\V{0},\V0}\right\rvert}\left\lvert\mean{\surr_{\jset\bigcup\{[\aindex,\bindex]\}}(\rvv)}-\mean{\surr_\jset(\rvv)}\right\rvert & & \derr^{\sigma^2} = \frac{1}{\left\lvert \surr_{\V{0},\V0}\right\rvert^2}\left\lvert\var{\surr_{\jset\bigcup\{[\aindex,\bindex]\}}(\rvv)}-\var{\surr_\jset(\rvv)}\right\rvert \nonumber\\ \gamma_\bl=\frac{1}{\dwork}\left(\kappa\derr^{\mu}+(1-\kappa)\derr^{\sigma^2}\right), & &\kappa\in[0,1]
\end{align}
The refinement indicator $\gamma_\bl$ attempts to maximize the reduction of the error in the mean and variance of the approximation whilst minimizing cost. Specifically the error indicator measures (via a convex combination) the contribution of the interpolant $\surr_\bl$ to the mean and variance of the MISC approximation $\surr_\jset$ which can be computed using \eqref{eq:general-misc}. 

We remark here that initialization challenges exist when using $\kappa=0$ to construct an adaptive MISC approximation. For any index $[\aindex,\bindex]$ with $\bindex=\V{0}$, the interpolants $\surr_{\aindex,\bindex}$ are built using only a single point and so are just constant functions. For these interpolants, the variance is zero and so not a useful indicator of the impact of refinement for initial $\aindex$ candidates.
Also note that the error indicators, $\derr^{\mu}$ and $\derr^{\sigma^2}$, are relative to the magnitude of the function and its square at the center of the sparse grid. This is important to ensure that mean or the variance components of the indicator dominate simply because they are on a different scales.

The index $\blb$ with the largest refinement indicator is refined using the function \texttt{REFINE} (step 6 of Algorithm \ref{alg:dim-adaptivity}) by adding all indices $\blb+\V{e}_k$, $k\in[\na+\nb]$ that satisfy the following downward-closed admissibility criterion
\begin{equation}
\label{eq:gsg_admissibility}
\blb+\mathbf{e}_k-\mathbf{e}_j\in\cL\text{ for } j\in[\na+\nb],\, l_k > 1
\end{equation}
The active set $\cA$ is then rebuilt by adding each index corresponding to the indices from \eqref{eq:gsg_admissibility}. 

\begin{figure}[ht]
\begin{center}
\includegraphics[width=.49\textwidth]{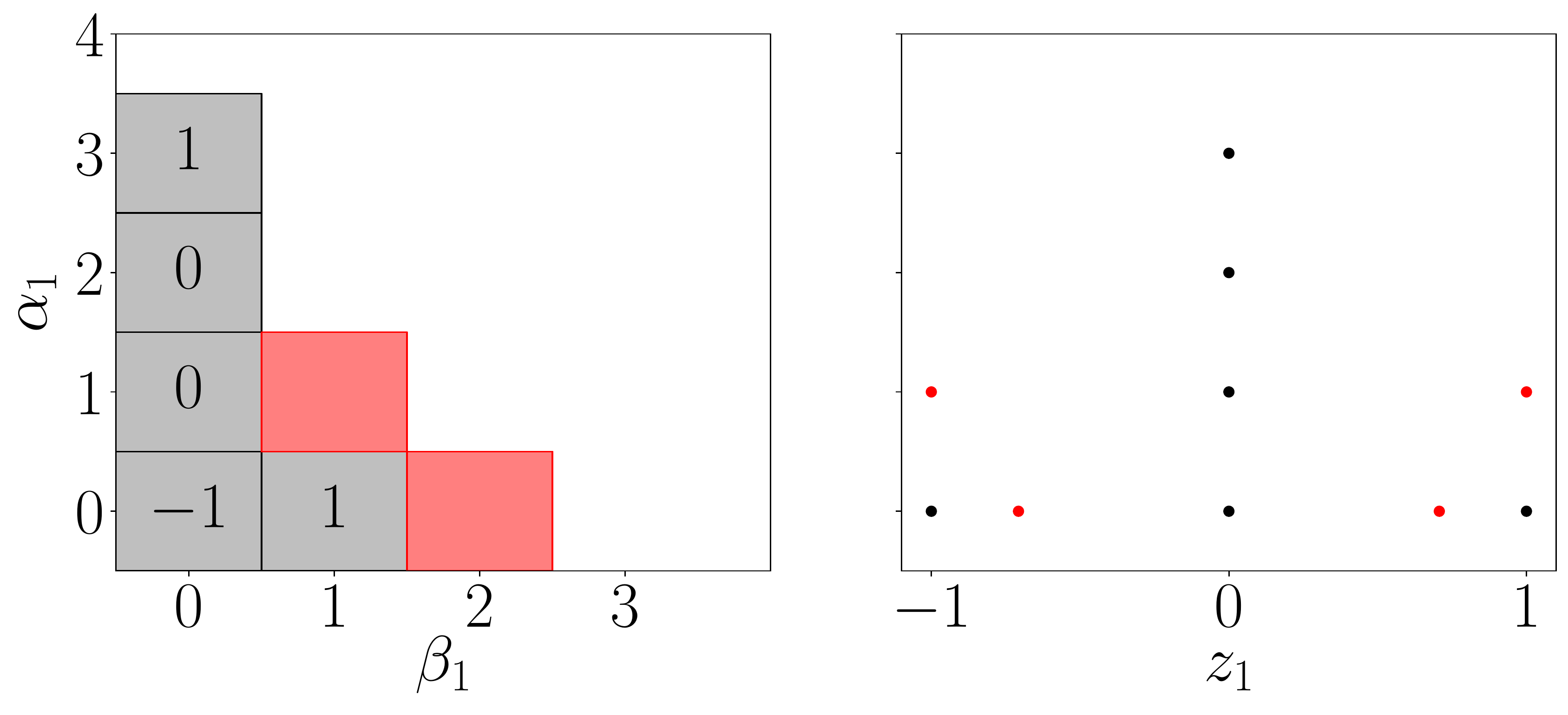}
\includegraphics[width=.49\textwidth]{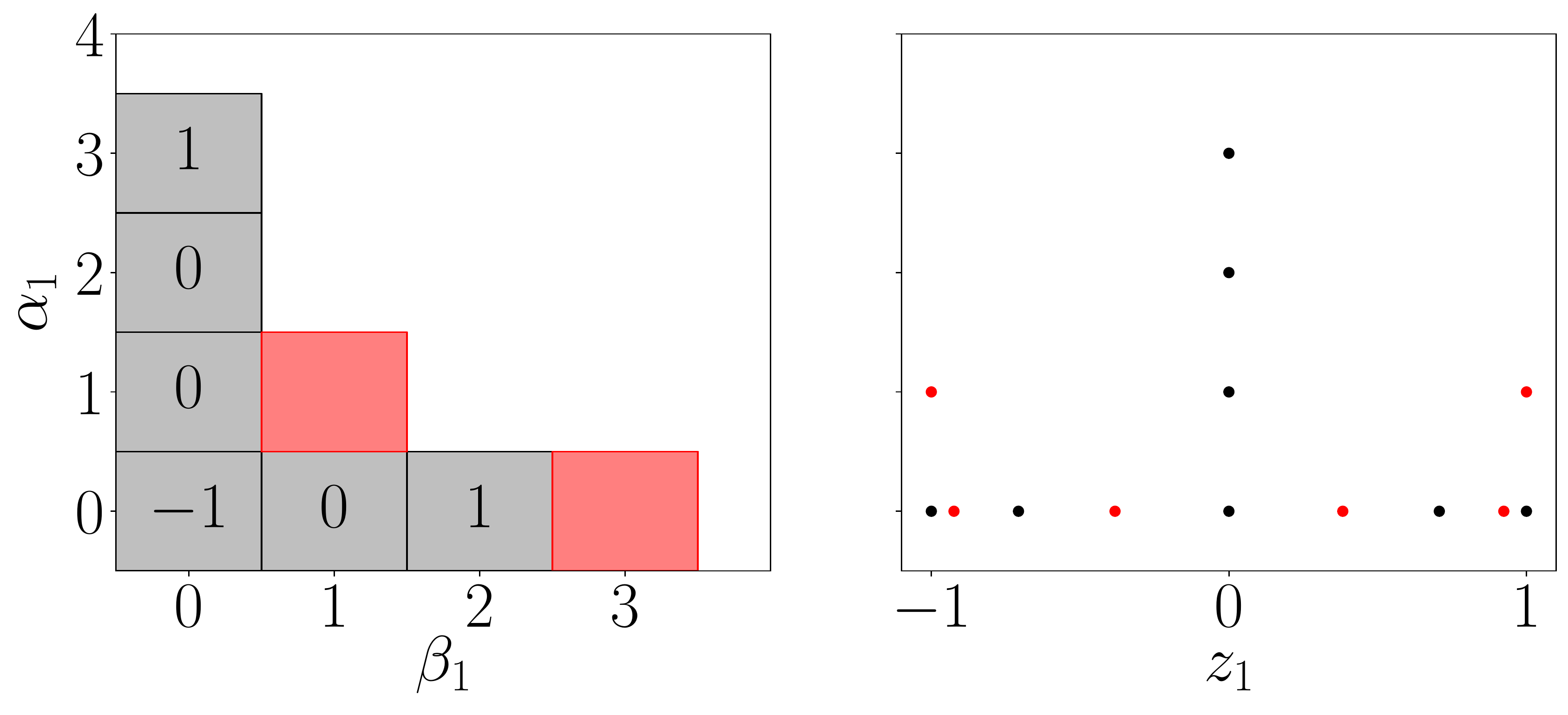}
\end{center}
\caption{Two steps of the adaptive sparse grid algorithm for $\na=\nb=1$. The gray and red boxes represent the sets $\jset$ and $\cA$, respectively. The numbers within the boxes are the sparse grid coefficients. The second and fourth plots represent the samples assigned to each level of physical discretization level $\ai_1$. The points corresponding to the tensor-product interpolations indexed by the set $\jset$ and the active set $\cA$ are given in black and red, respectively.}
\label{fig:adaptive-misc}
\end{figure}

Figure \ref{fig:adaptive-misc} depicts two steps of the adaptive sparse grid algorithm for $\na=\nb=1$. The first and third plots depict the index sets $\jset$ and $\cA$ represented by the gray and red respectively. The numbers within the boxes are the sparse grid coefficients. The second and fourth plots represent the samples assigned to each level of physical discretization level $\ai_1$, $l_{\ai_1}=0,\ldots,3$. The points corresponding to the tensor-product interpolants indexed by the set $\jset$ and the active set $\cA$ are given in black and red respectively. The new work required to refine an index is given by the cost of running the simulation at all the red points $\diffgrid_{\aindex,\bindex}$.  In a previous step, the index $(\ai_1,\bi_1) = (0,1)$ was chosen for refinement. Both forward neighbors $\V{k} = \blb+\V{e}_n$, $n\in[\na+\nb]$ are added to the active set $\cA$ (red boxes in first figure) because these indices satisfy the admissibility criteria. This requires evaluating both the functions $\surr_{0}$ and $\surr_{1}$ at 2 new points. In the next step, the index $(\ai_1,\bi_1) = (0,2)$ is selected for refinement. In this step, only one new index  $(\ai_1,\bi_1)=(0,2)+(0,1)=(0,3)$ is added to the active set (third figure) because the other forward neighbor $(\ai_1,\bi_1)=(0,2)+(1,0)=(1,2)$ does not satisfy the admissibility criteria. Adding the index $(2,0)$ to $\cA$ requires evaluating $\surr_{0}$ at 4 new points.

\rev{\subsection{Multi-level Collocation}
\label{sec:multi-level-collocation}

Multi-level collocation (MLC) ~\cite{Teckentrup_JWG_SIAMUQ_2015} is a specialized instance of multi-index collocation. MLC assumes that there is only one discretization index which controls the model accuracy and cost, i.e. $\aindex=(\ai_1)$. We can easily apply the AMISC algorithm to the multi-level setting by simply defining the index set $\mathcal{M}_\aindex$ correctly. For the advection-diffusion considered in this paper we make the reasonable choice that $\mathcal{M}_\aindex=\{\otimes_{i=1}^3\{1\},\ldots,\otimes_{i=1}^3\{k\}\}$, $k\in\{1,\ldots,5\}$. This multi-index model set has a one-to-one mapping with the one-dimensional multi-level hierarchy $\hat{\mathcal{M}}_\aindex=\{1,2,3,4,5\}$. We do not consider multi-level methods for the nozzle problem as it is unclear how to define a reasonable one-dimensional hierarchy for that model.

The original multi-level algorithm \cite{Teckentrup_JWG_SIAMUQ_2015} allocate samples to models of different model discretizations using a priori theoretical error estimates. Consequently it is challenging to apply this approach to practical problems where such theoretical estimates are unavailable. A partially adaptive a posteriori MLC strategy was proposed in \cite{Farcas_BNU_SGA_2018}. However, when using this approach, the resolution of the approximation built over the stochastic space is a priori tied to the resolution of the numerical solver used to solve the PDE. The adaptive multi-index collocation presented here can easily be applied in the multi-level setting without the restrictions of the aforementioned approaches. The computational cost and accuracy of multi-level and multi-index sparse grid collocation are compared in Section~\ref{sec:results-advection-diffusion}.

Finally it is important to note that multi-level surrogate methods have been developed which do not use sparse grids to build multivariate approximations over the stochastic space. Gaussian processes were used in \cite{Perdikaris_PRSLA_2016,Kennedy_O_B_2000,LeGratiet_G_IJUQ_2014} and low-rank reduced order models were used in \cite{Narayan_GX_SISC_2013}. These methods are not considered further here since, similar to multi-level sparse grid collocation, if a model does indeed have multiple discretization hyper-parameters it can be unclear how to define a suitable one-dimensional model hierarchy.
}

\section{Multi-index sensitivity analysis}
\label{sec-5}
\label{sec:sensitivity-analysis}
Quantifying the sensitivity of a model output $f$ to the model parameters $\rVv$ can be an important component of assessing uncertainty in model predictions. For example, parameter sensitivities can be used to prioritize data collection. Using data or additional knowledge to reduce uncertainty in sensitive model inputs can significantly decrease prediction uncertainty. This section demonstrates how to use MISC for global sensitivity analysis, specifically for computing Sobol indices~\cite{Sobol_MMCE_1993} . 

\rev{Sobol indices are frequently used to estimate the sensitivity of a function to single, or combinations of, input parameters. Sobol indices can be estimated numerically using Monte Carlo methods \cite{Saltelli_book_2000}; however, they can be computed analytically from a sparse grid \rev{\cite{Buzzard_RESS_2012,Formaggi_GILPRST_CG_2013}}. For any tensor-product probability density function $\pbwt$, there exists a transformation of the Lagrange polynomial basis into a orthonormal polynomial basis \cite{Buzzard_CBJ_2013}, also known as a polynomial chaos expansion (PCE). Applying such a transformation facilitates simple and analytical computation of Sobol indices from the coefficients of the PCE \cite{Sudret_RESS_2008}. The following sections describe the basic properties of a PCE and how to transform an MISC approximation into a PCE for the purpose of computing Sobol indices.}

\subsection{Polynomial chaos expansions}
\label{sec-5-1-1}
Polynomial chaos expansions (PCE) \cite{Ghanem_book_1991,Jakeman_FNEP_CMAME_2019,Xiu_K_SISC_2002} represent the model output \(\ffem(\rvv)\) as an expansion of orthonormal polynomials
\begin{align}
\label{eq:pce-multi-index}
\ffem(\rvv) &\approx \surr_{\aindex,\rev{\Lambda}}(\rvv) = \sum_{\lambda\in\Lambda}\rev{\eta}_{\lambda}\phi_{\lambda}(\rvv), & |\Lambda| &= N.
\end{align}
The basis functions $\phi_\lambda$ are typically constructed to be orthonormal with respect to the density \(\pbwt\), that is 
\[
  (\phi_{\lambda^{(i)}}(\rvv),\phi_{\lambda^{(j)}}(\rvv))_{L^2_{\pbwt}(\Omega)} \coloneqq \int_{\rvvsupp} \phi_{\lambda^{(i)}}(\rvv)\phi_{\lambda^{(j)}}(\rvv) \dx{\pbwt}(\rvv)= \delta_{i,j},
\]
where $\delta_{i,j}$ is the Kronecker delta function and $\lambda^{(i)},\lambda^{(j)}$ are two different multivariate indices where \(\lambda=(\lambda_1\ldots,\lambda_\nv)\in\mathbb{N}_0^\nv\). 
When the components of \(\rvv\) are independent we can express the basis functions \(\phi\) as tensor products of univariate orthonormal polynomials. That is
\begin{align}
\phi_\lambda(\rvv)=\prod_{i=1}^\nv \phi^i_{\lambda_i}(\rv_i), & &   \int_{\Omega_i} \phi^i_{j}(z_i) \phi^i_{k}(z_i) \pbwt_i(z_i) \dx{z_i} &= \delta_{j,k}, & j, k &\geq 0, & \deg \phi^i_j &= j.
\end{align}

The mean and variance of the PCE can be computed directly from the 
expansion coefficients.
\begin{align}
\label{eq:gpc_mean_var}
\mean{\surr_{\aindex,\rev{\Lambda}}}=\rev{\eta}_\V{0}, & & \var{\surr_{\aindex,\rev{\Lambda}}}=\sum_{\substack{\V{\lambda}\in\Lambda\\ \V{\lambda}\neq \V{0}}} \rev{\eta}^2_\V{\lambda}
\end{align}
The Sobol indices can also be computed analytically from the expansion using the identity
\begin{align}
\label{eq:pce-sobol-indices}
\var{\surr_\V{u}}=\sum_{\V{\lambda}\in\Lambda_\V{u}}\rev{\eta}_\V{\lambda}^2 & & \Lambda_\V{u}=\{\lambda \mid \lambda_{i}>0,\; i\in\V{u},\; \lambda_j=0,\; j\notin \V{u}\}
\end{align}

\subsection{Transforming a sparse grid into a polynomial chaos expansion}
\label{sec-5-1-2}
Let $\C{T}$ denote the transformation of a Lagrange polynomial interpolant $\surr_{\aindex,\bindex}$ to an orthonormal polynomial interpolant $\surr_{\aindex,\rev{\Lambda}}$. From \eqref{eq:general-misc}, we have
$\C{T}[\misc(\rvv)]
=\sum_{[\aindex,\bindex]\in\jset}c_{\aindex,\bindex}\,\C{T}[\surr_{\aindex,\bindex}(\rvv)]$, where $c_\bindex$ are the Smolyak coefficients. Thus, to transform the sparse grid, we only need the ability to transform each tensor-product interpolation individually.
To apply the transformation of the tensor-product interpolants, we must only compute transformations of one-dimensional basis functions. In one dimension, the polynomial interpolant of a set of points is unique. Consequently, we can express the $j$-th univariate Lagrange basis function in the $i$-th variable dimension as 
\begin{align}
l_{i,j}(\rv_i)=\sum_{k=1}^{m_{\bi_i}}\nu_k\phi^i_k(\rv_i) && \nu_k=\int_{\rvvsupp_k} l_{i,j}(\rv_i)\phi^i_k(\rv_i)\dx{\pbwt(\rv_i)} && j=1,\ldots,m_{\bi_i}, i=1,\ldots,\nv
\end{align}
We can compute the coefficients of the orthonormal basis using Gaussian quadrature. Note this step does not require any evaluations of the model $\ffem$ and can be precomputed offline and stored for future use. Given the orthonormal representation of the univariate Lagrange polynomials, the tensor-product interpolant can be expressed as 
\begin{align}
\label{eq:multivariate_lagrange-pce-transformation}
\surr_{\aindex,\bindex}(\rvv) = \sum_{\V{j}\le\bindex}
\ffem(\rvv^{(\V{j})})\prod_{i=1}^dl_{i,j_i}(\rv_i)=
\sum_{\V{j}\le\bindex}
\ffem(\rvv^{(\V{j})})\prod_{i=1}^d\left(\sum_{k=1}^{m_{\bi_i}}\nu_k\phi^i_k(\rv_i)\right)=\sum_{\V{\lambda}\in\Lambda_\bindex}\rev{\eta}_\lambda \phi_\V{\lambda}(\rvv),
\end{align}
where $\Lambda_\bindex=\{\V{\lambda}\mid \V{\lambda}\le\bindex\}$. 
This expression can be now used with \eqref{eq:pce-sobol-indices} to compute Sobol indices analytically and efficiently from the sparse grid.

Note that the algorithm described here is different from the algorithm\rev{s} proposed in \rev{\cite{Buzzard_RESS_2012,Formaggi_GILPRST_CG_2013}. These approaches} solve a multivariate interpolation problem to compute the coefficients $\gamma$ of the tensor-product PCE basis in \eqref{eq:multivariate_lagrange-pce-transformation}. The dominant cost of that approach is the inversion of the multivariate Vandermonde-like matrix evaluated at each of the tensor product points, which requires $O(M^3(\beta))$ operations. By focusing on univariate transformations of the Lagrange basis, we can compute the coefficients of the multivariate PCE via an outer product of the univariate orthonormal basis coefficients at a cost of $O(M^2(\beta))$. This reduction in the cost of computing the PCE transformation is practically significant, but the complexity of the sparse grid transformation still grows linearly with the number of points in the sparse grid. Note that the cost of either transformation algorithm is typically negligible to the cost of evaluating the models $\ffem$.

\begin{remark}
  In this paper, we use sparse grids to form a MISC approximation. Sparse grids naturally combine the refinement of the discretization parameters and the refinement of the interpolants. Noting that the MISC approximation can be expressed in the general form $\misc(\rvv) = \sum_{[\aindex,\bindex]\in\jset}c_{\aindex,\bindex}\,\surr_{\aindex,\bindex}(\rvv)$, equation \eqref{eq:multivariate_lagrange-pce-transformation} suggest that sparse grids are not the only means to build a MISC approximation. Specifically, one could draw from a broad class of surrogate models to construct the approximations $\surr_{\aindex,\bindex}$. One could use PCE constructed using interpolation, least-squares or $\ell_1$ minimization, but one could also use Gaussian processes, low-rank approximations, or even neural networks. Given a method to build $\surr_{\aindex,\bindex}$, the adaptive algorithm described in Section \ref{sec:adaptivity} can be used to guide the selection of the indices $[\aindex,\bindex]$ so long as a sampling method, which can be controlled by the stochastic discretization parameter $\bindex$, is available. It is straightforward to construct such sampling methods; for example, one approach is to let $\bindex$ be a scalar controlling the number of Monte-Carlo or Latin-hypercube samples used to build a regression-based PCE or Gaussian process. Alternatively, one can use sampling schemes such as multivariate Leja sequences for polynomial interpolation, which are more amenable to dimension-adaptivity~\cite{Jakeman_FNEP_CMAME_2019}.
\end{remark}

\section{Results}
\label{sec-6}
\label{sec:results}
In this section, we explore the efficacy and properties of our proposed approach, first using the advection-diffusion model problem and then using the engineering model of a jet engine nozzle.

To measure the performance of an approximation, we will use the \(\ell^\infty\) error on a set of test nodes.  We generate a set of \(Q\) random samples \(\{\samp{j}\}_{j=1}^Q\subset\rvvsupp\) drawn from the density \(\pbwt\) of the uncertain variables $\rVv$. The relative error is computed as 
\[
  \label{eq:error}
\norm{f-\surr_\jset}{L^\infty(\rvvsupp)} = \frac{1}{\max_{j\in[Q]} f(\samp{j})-\min_{j\in[Q]} f(\samp{j})}\argmax_{j\in[Q]} \abs{f(\samp{j})-\surr_\jset(\samp{j})}
\]
where $f$ is the exact function and $\surr_\jset$ is the interpolative approximation. We set $Q=1000$ and $Q=500$ for the advection diffusion and nozzle models respectively.

We use the  \(\ell^\infty\) norm to quantify error as we are interested in probabilities of failure, and convergence almost surely is required to guarantee convergence of probability density functions (PDFs) of the QoI \cite{Butler_JW_SISC_2018}. We use the range of the values of the reference function $f$ over the validation samples to normalize the error, enabling us to highlight how our adaptive approximations balance error among multiple QoI.

\subsection{Advection diffusion model}
\label{sec:results-advection-diffusion}
In this section, we will use the adaptive MISC algorithm to approximate
the quantity of interest \eqref{eq:advec-diff-qoi} obtained from the advection diffusion model \eqref{eq:advection-diffusion}. For this model, an efficient and accurate approximation must balance spatial refinement, temporal refinement, and the stochastic interpolation error. The details of the physical discretization are provided in Section \ref{sec:adv-diff-phys-discretization}.

\begin{figure}
\centering
\includegraphics[width=0.49\textwidth]{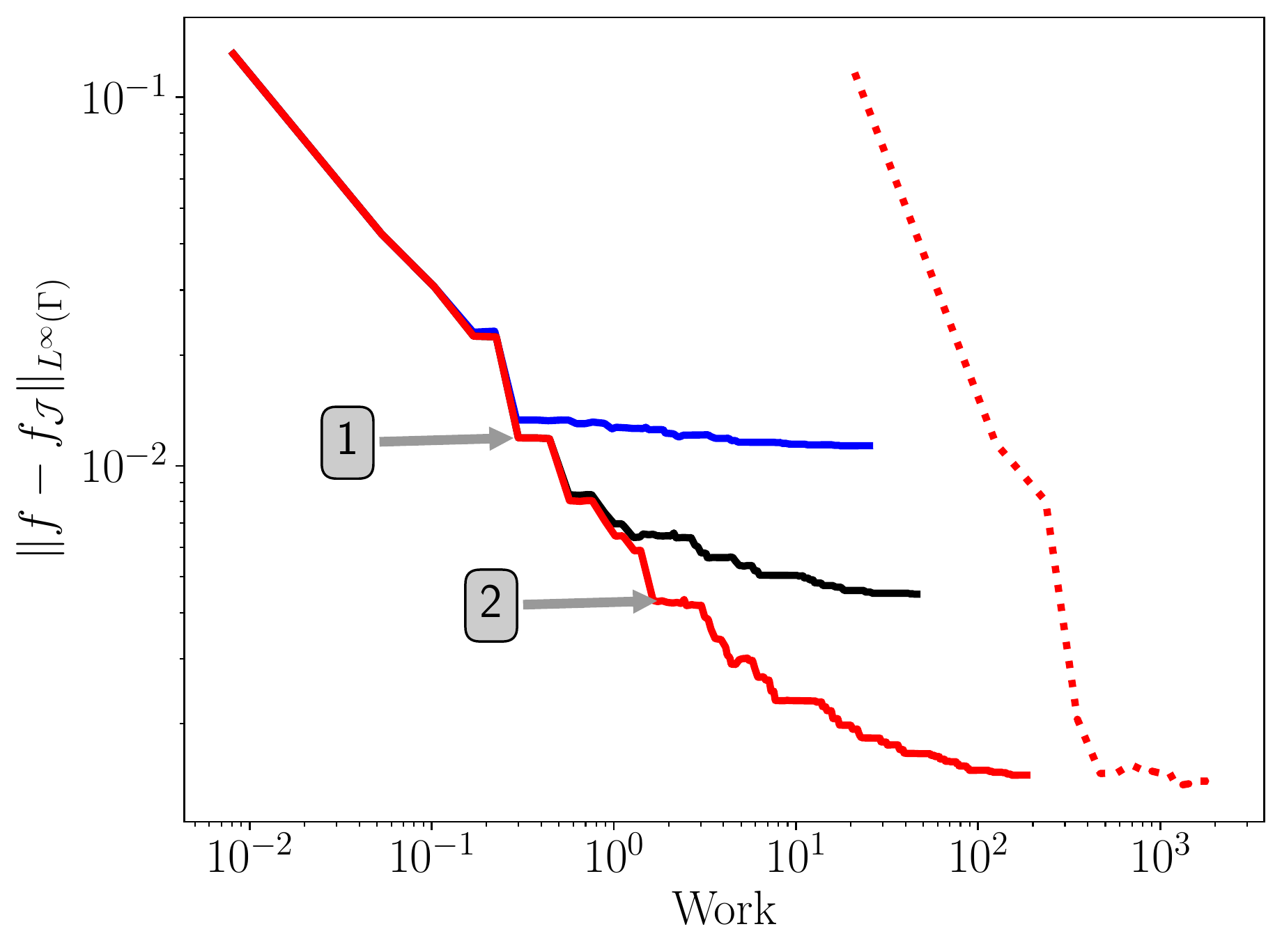}
\includegraphics[width=0.49\textwidth]{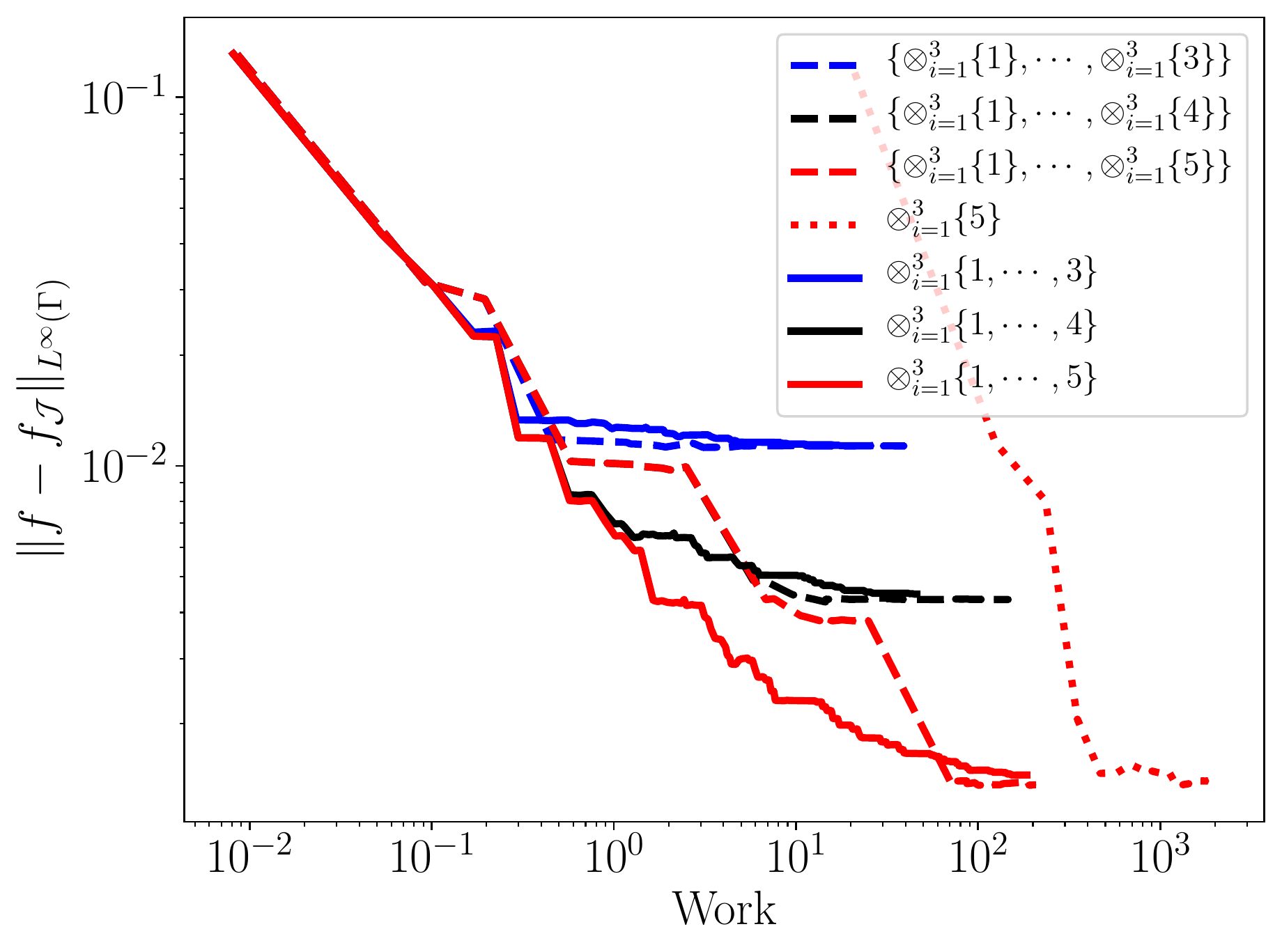}
\caption{\rev{(Left) Comparison of the convergence of single-fidelity and multi-fidelity approximations of the advection-diffusion model. (Right) Comparison of AMISC with adaptive multilevel approximation. The legend labels denote the index sets $\mathcal{M}_\aindex$ defining the ensemble of models considered by the different algorithms.  $\mathcal{M}_\aindex=\otimes_{i=1}^3\{5\}$ denotes the single fidelity approximation using evaluations of only the highest fidelity model $\surr_{(5,5,5)}$, $\mathcal{M}_\aindex=\otimes_{i=1}^3\{1,\ldots,k\}$ and $\mathcal{M}_\aindex=\{\otimes_{i=1}^3\{1\},\ldots,\otimes_{i=1}^3\{k\}\}$ respectively denote multi-index and multi-level approximation where the highest fidelity model evaluated is $\surr_{(k,k,k)}$, $k\in\{3,4,5\}$.}}
\label{fig:adv-diff-convergence}
\end{figure}
\begin{figure}[ht]
\begin{center}
\includegraphics[width=\textwidth]{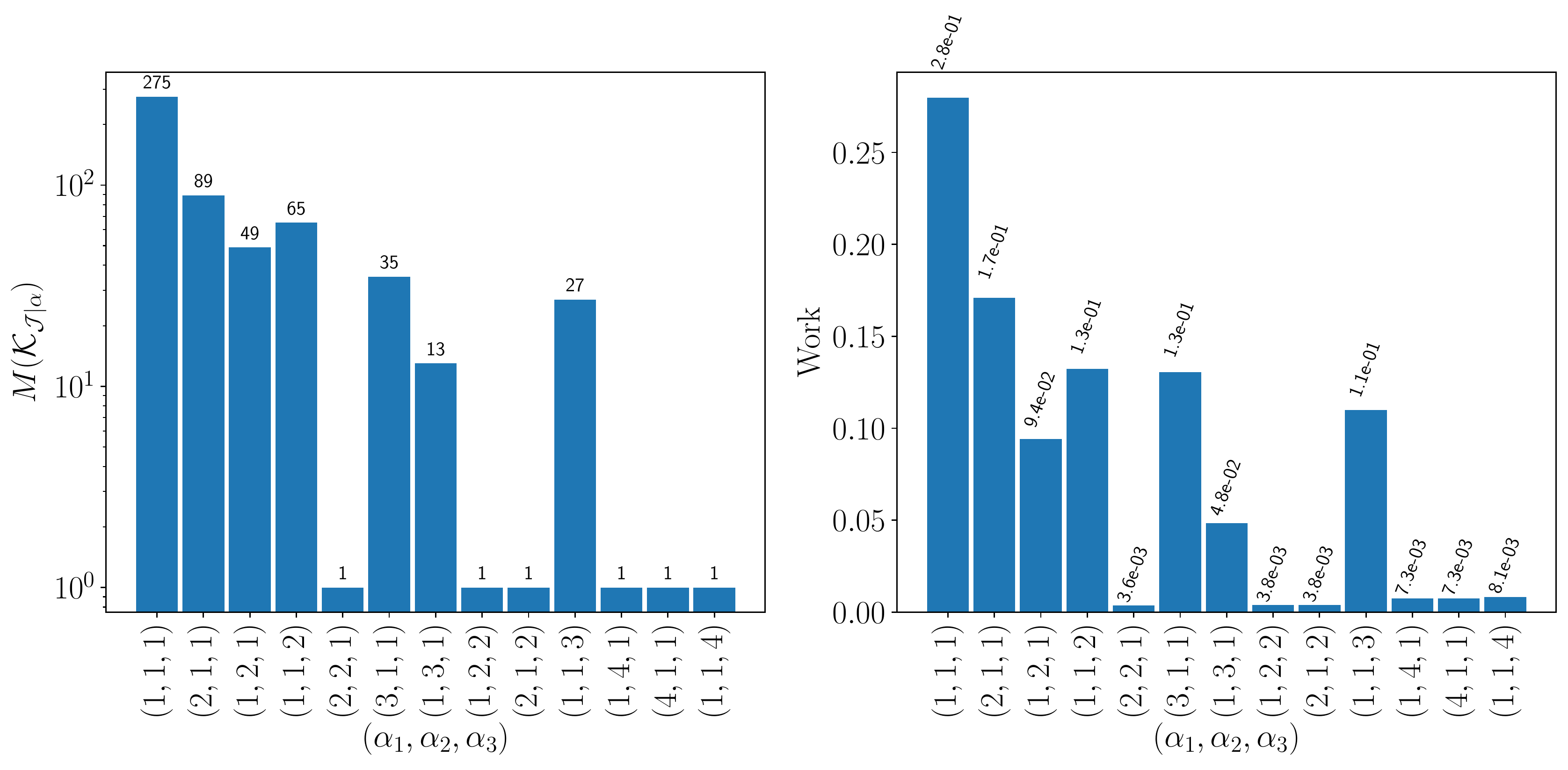}\\
\end{center}
\caption{The number of simulations (left) and fraction of total work (right) assigned to each approximate advection diffusion model at the point labeled 1 in Figure \ref{fig:adv-diff-convergence}. \rev{Here $M_\kset=\mathrm{card}\left(\grid_{\kset}\right)$ is the number of samples, of the random variables, at which $\surr_\aindex$ is evaluated. The number on top of the bars in the left plot are the exact values of $M_\kset$ and the numbers on the bars on the right plot are the fraction of total work
assigned to each model indexed by $\aindex$. Finally the labels of the horizontal axis represent the indices $\aindex$ to
which AMISC has assigned at least one evaluation.}}
\label{fig:adv-diff-cost-profile-1}
\end{figure}
\begin{figure}[ht]
\begin{center}
\includegraphics[width=\textwidth]{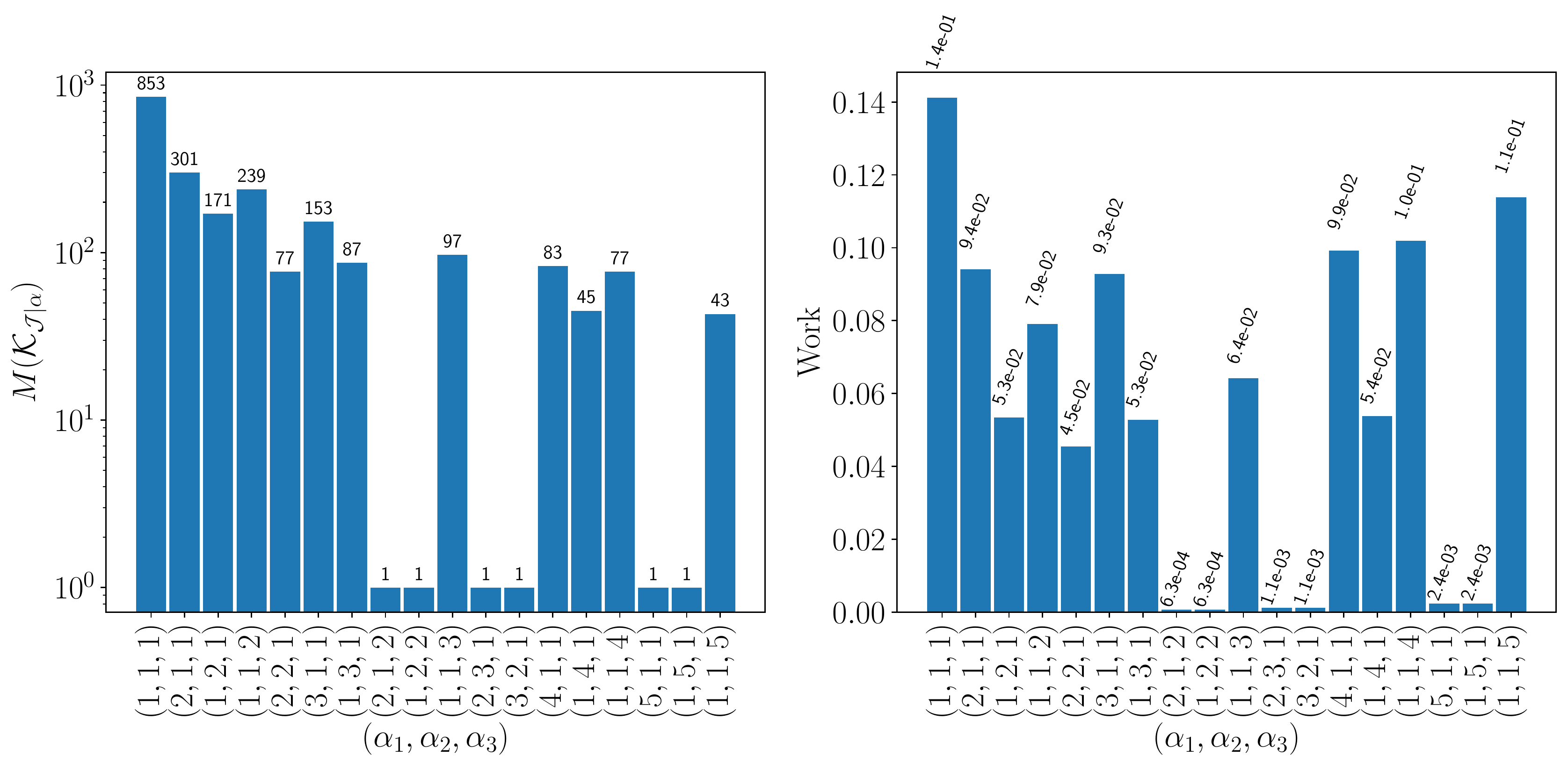}
\end{center}
\caption{The number of simulations (left) and fraction of total work (right) assigned to each approximate advection diffusion model at the point labeled 2 in Figure \ref{fig:adv-diff-convergence}. \rev{Here $M_\kset=\mathrm{card}\left(\grid_{\kset}\right)$ is the number of samples, of the random variables, at which $\surr_\aindex$ is evaluated. The number on top of the bars in the left plot are the exact values of $M_\kset$ and the numbers on the bars on the right plot are the fraction of total work
assigned to each model indexed by $\aindex$. Finally the labels of the horizontal axis represent the indices $\aindex$ to
which AMISC has assigned at least one evaluation.}}
\label{fig:adv-diff-cost-profile-2}
\end{figure}

In Figure \ref{fig:adv-diff-convergence}, we compare the adaptive MISC algorithm (AMISC) with a single level adaptive sparse grid (see Section \ref{sec:sparse-grid}). Specifically, we plot the error in each approximation as the total work is increased. Work is measured relative to the cost of performing a single simulation using the approximate model $\surr_{5,5,5}$ with the highest fidelity discretization parameters. \rev{For example, if $\mathrm{\tau}_\aindex$ is the wall time needed to simulate the model $f_{\aindex}$ at one realization of the parameters $\rvv$, then the work needed to run this model is $W_\aindex=\tau_\aindex/\tau_{5,5,5}\le 1$.
For this particular advection-diffusion problem, we do not know the true QoI exactly, so we use the model with $\aindex=(6,6,6)$ as the truth when computing errors.
Each time the error is measured, all points associated with indices in both the active index set $\cA$ and the current index set $\iset$ are used and thus contribute to the value of work on the horizontal axis. The notation in the legend is used to denote the index sets $\mathcal{M}_\aindex$ defining the ensemble of models considered by the different algorithms.  $\mathcal{M}_\aindex=\otimes_{i=1}^3\{5\}$ denotes the single fidelity approximation using evaluations of only the highest fidelity model $\surr_{(5,5,5)}$, $\mathcal{M}_\aindex=\otimes_{i=1}^3\{1,\ldots,k\}$ and $\mathcal{M}_\aindex=\{\otimes_{i=1}^3\{1\},\ldots,\otimes_{i=1}^3\{k\}\}$ respectively denote multi-index and multi-level approximation where the highest fidelity model evaluated is $\surr_{(k,k,k)}$, $k\in\{3,4,5\}$. The gray boxes refer to two iterations during the
evolution of the MISC algorithm. The point labeled 1 is taken when 500 model evaluations (of any fidelity)
have been run and label 2 after 2000 model evaluations. The cost profiles of the MISC algorithm at these two iterations are depicted in Figures \ref{fig:adv-diff-cost-profile-1} and \ref{fig:adv-diff-cost-profile-2}.}

The error in both single-level sparse grid and AMISC approximations eventually saturate when the physical discretization error of the highest fidelity approximate model is reached. As the maximum value of $\ai_1$ is increased, the error in the approximations saturates at a smaller error. The AMISC approximations, however, achieve the same level of error as the single fidelity approach at only a fraction of the work. The results in Figure \ref{fig:adv-diff-convergence} are generated by setting $\kappa=1$ in \eqref{eq:dim-surplus-indicator}. Despite the fact we only adapt to reduce the error in the mean, we still can reduce the $\ell^\infty$ error in the AMISC and single-fidelity approximations.

The AMISC algorithm reduces the total computational cost by over two orders of magnitude\rev{, except in regions where the physical discretization error of the highest fidelity model begins to dominate}.\footnote{\rev{If the cost profile of a set of model resolutions increases more rapidly, then the computational benefit of AMISC also increases. For example, if the conservative cost model for solving the advection-diffusion problem presented in Section~\ref{sec:adv-diff-phys-discretization} is replaced with a cost model that grows more quickly with the mesh and temporal resolution, then the benefits of the multi-fidelity methods will increase relative to their single fidelity counterparts.}}
 This reduction is achieved by sampling the lower-level model approximations more than the higher-level models, and only sampling the high-fidelity models when the physical discretization error starts to dominate the stochastic error. Figures \ref{fig:adv-diff-cost-profile-1} and \ref{fig:adv-diff-cost-profile-2} plot the total number of simulations the MISC algorithm assigns to each approximate model for two points in the algorithms evolution (labeled 1 and 2 in Figure \ref{fig:adv-diff-convergence}). \rev{When the stochastic error is approximately equal to the saturation point of the blue curve, most simulations are allocated to the low-fidelity models which only involve models with indices $\ai_i\le3,i=1,\ldots,3$ (see Figure \ref{fig:adv-diff-cost-profile-1}). However when the stochastic error is below the black curve, numerous samples are allocated to the higher-fidelity models with indices $\ai_i>3$ (see Figure \ref{fig:adv-diff-cost-profile-2}).} Although most evaluations are assigned to the lowest-level models, \rev{the highest fidelity models contribute significantly to the total work (Figure \ref{fig:adv-diff-cost-profile-2} right) despite being evaluated much less frequently}. This is the behavior we desire. Specifically, the high-fidelity models are only evaluated when needed to maintain deterministic prediction accuracy.

\rev{The right plot of Figure~\ref{fig:adv-diff-convergence} compares multi-index collocation with multi-level collocation using the multi-level hierarchy defined in Section~\ref{sec:multi-level-collocation}. Multi-index collocation obtains significant computational savings when high-levels of accuracy are required. The computational savings increase as the number of levels utilized by the multi-level method increase. Furthermore, if the cost differential across a set of model resolutions becomes more significant, the computational benefit of MISC compared to MLC will increase.

The results generated using MLC were found using, what turned out to be, a reasonable one-dimensional model hierarchy. It is easy to construct other hierarchies which significantly degrade the performance of MLC relative to AMISC. One major benefit of AMISC is that when a model does have multiple discretization hyper-parameters, one does not need to decide how to form a one-dimensional model hierarchy.

As expected, the multi-level and multi-index approximation saturate at the same points when fixing the highest fidelity model used in the approximations consistently. The multi-index method does take longer to reach the exact saturation point, however the error gets close to the saturation point much more quickly than
the multi-level method. The small inefficiencies of the multi-index near the saturation points are caused by the downward-closed index set constraint for cases where final reductions in deterministic discretization error could only be addressed using the most resolved model, resulting in some loss in relative benefits for the tails of the convergence histories. This effect can be avoided by increasing the resolution of the highest fidelity model. Doing so will not effect the computational cost of the AMISC algorithm unless the evaluations of the higher fidelity model is needed to reduce the deterministic prediction error. Future work will leverage this fact to explore extensions to AMISC approaches that do not require the highest fidelity model to be specified a priori.}

\subsection{Aero-thermal-structural analysis of a jet engine nozzle}
\label{sec-6-2}
In this section, we will use the adaptive MISC algorithm to approximate the four quantities of interest obtained from the nozzle model, and then use this approximation to quantify sensitivities and uncertainty in these QoI. The four model QoI we consider are mass, thrust, load-layer temperature failure ratio, and thermal-layer failure criteria. For both failure-based QoI, a value greater than one indicates failure, and both are computed using the modified P-norm function (PN-function), which aggregates pointwise values into a global quantity:
\begin{equation}
\label{eq:pn_function}
g(c_i) = \left( \frac{1}{N} \sum_i^N c_i^p \right)^{\frac{1}{p}}
\end{equation}
As $p$ increases, the PN-function approaches the value $\textrm{max}_i(c_i)$, where $c_i$ denote local stress or temperature values.  Here, we set $p=10$.

\subsubsection{AMISC approximation}
\label{sec-6-2-1}
In this section, we detail the performance of the AMISC method for approximating the aforementioned QoI of the nozzle model. When considering our nozzle model, an efficient and accurate approximation must balance CFD and structural model errors. The relative computational cost of each nozzle model approximation is given in Table \ref{tab:nozzle-cost}, where cost is reported in wall time using 4 CPUs for the CFD model and 1 CPU for the structural model. 

For this problem, we are interested in approximating four QoI, so we must adjust slightly the refinement indicator introduced in Section \ref{sec:adaptivity}. Specifically, we use a worst case error indicator that refines the multi-index set in the AMISC approximation according to the largest error across the set of QoI. Letting $\gamma_\bl^i$ denote the error indicator \eqref{eq:dim-surplus-indicator} \rev{for the $i$th QoI}, we use
\begin{align}
\gamma_\bl=\max_{i\in[4]} \gamma_\bl^i
\end{align}
across the set of $4$ QoI.  We select $\kappa=1/2$ in order to balance error reduction for the mean and variance of the approximation. The latter is consistent with our desire to use variance-based decomposition for global sensitivity analysis. The former ensures that the refinement indicator is not zero for constant functions, including approximations for $\bindex=\V{0}$ (see the comment in Section \ref{sec:adaptivity}). Also note that the use of a relative error indicator is important to ensure that the thrust does not dominate refinement since its scale is orders of magnitude larger than the other QoI (Figure \ref{fig:qoi-pdfs} shows the differing scales for the four QoI).

\begin{figure}[htb]
\begin{center}
\includegraphics[width=0.4\textwidth]{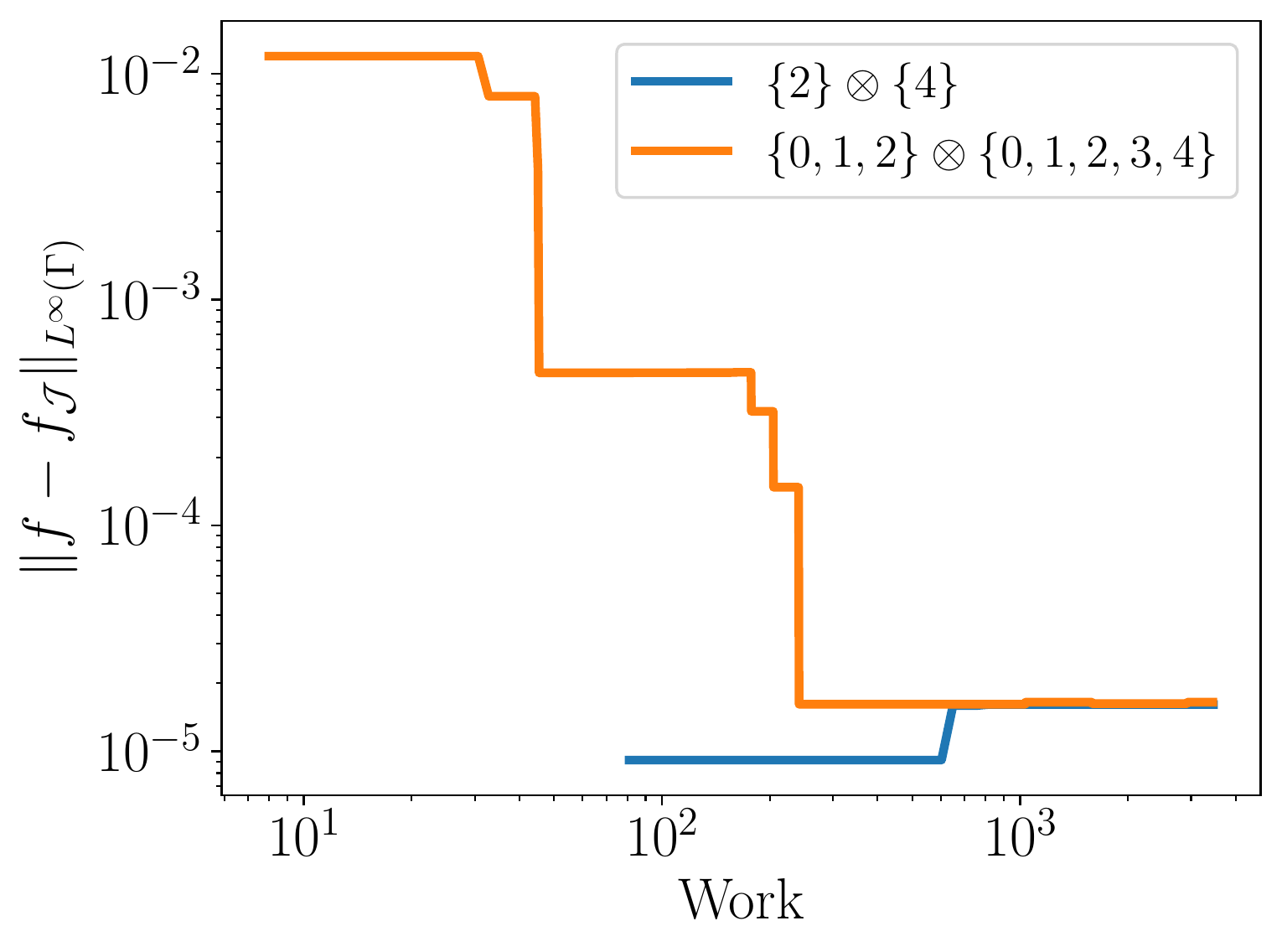}
\includegraphics[width=0.4\textwidth]{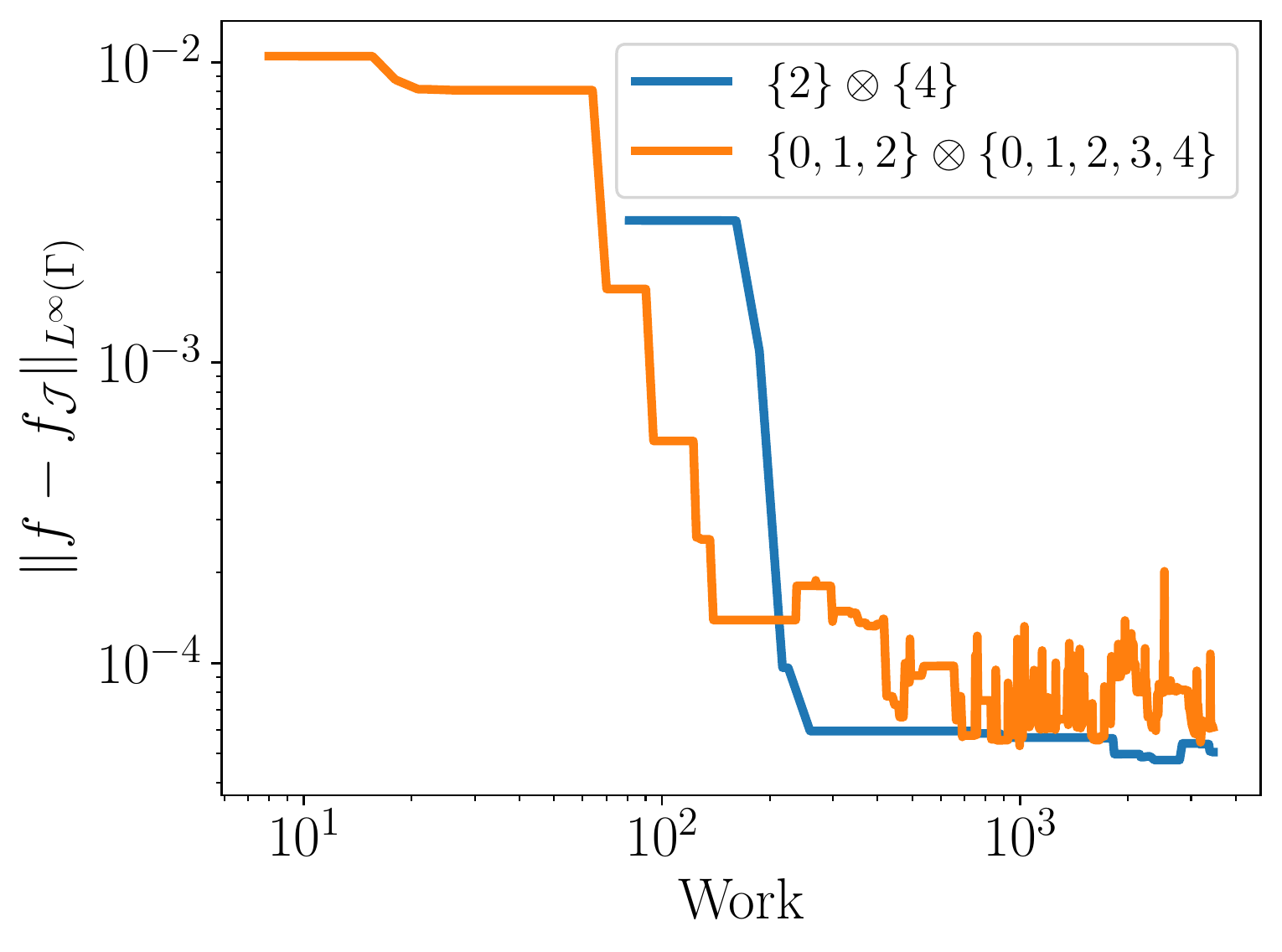}
\includegraphics[width=0.4\textwidth]{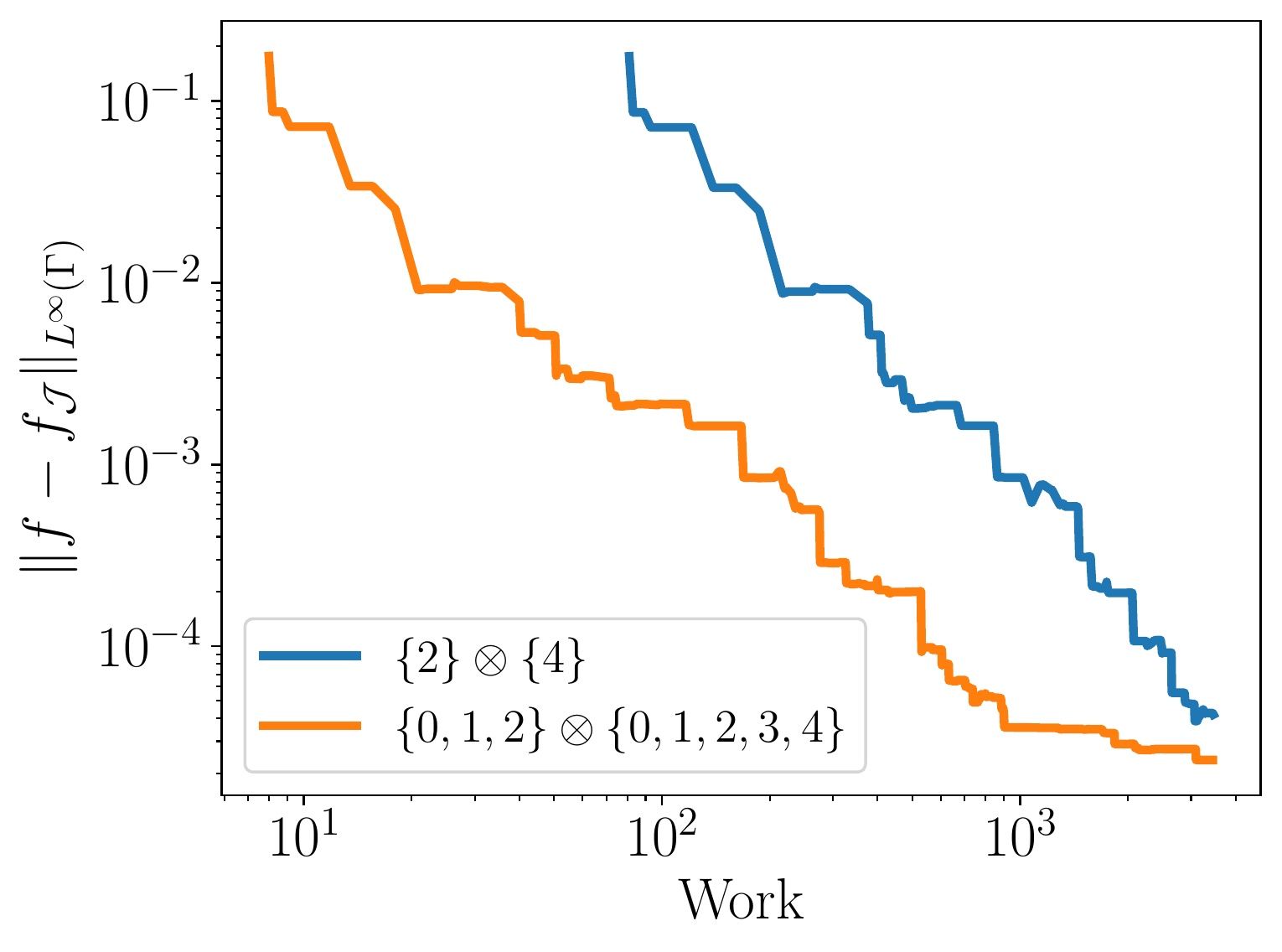}
\includegraphics[width=0.4\textwidth]{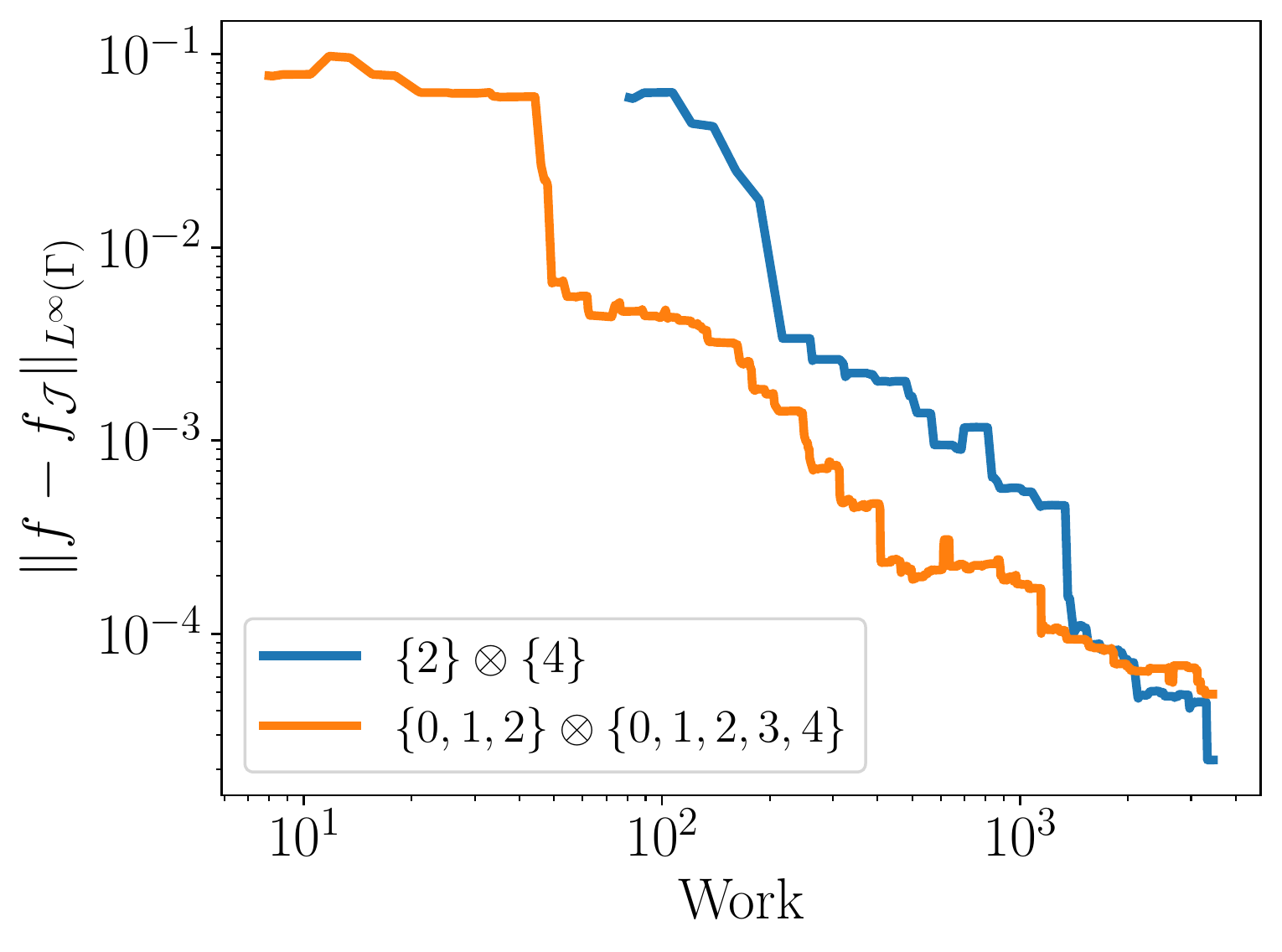}
\end{center}
\caption{Comparison of the convergence of single-fidelity and multi-fidelity approximations of (clockwise from top left) mass, thrust, temperature layer failure and load layer failure. \rev{The legend labels denote the index sets $\mathcal{M}_\aindex$ defining the ensemble of models considered by the different algorithms.  $\mathcal{M}_\aindex=\{2\}\otimes\{4\}$ denotes the single fidelity approximation using evaluations of only the highest fidelity model $\surr_{(2,4)}$, and $\mathcal{M}_\aindex=\{0,1,2\}\otimes\{0,1,2,3,4\}$ denotes the multi-index approximation where the highest fidelity model evaluated is $\surr_{(2,4)}$.}}
\label{fig:nozzle-error-convergence}
\end{figure}

In Figure \ref{fig:nozzle-error-convergence}, we compare the adaptive MISC algorithm (AMISC) with a single-level adaptive sparse grid (see Section \ref{sec:sparse-grid}). Work is measured relative to the cost of simulation using the highest fidelity discretization parameters (model $\surr_{3,5}$). We terminate refinement of each approximation after the cost of 3400 high-fidelity samples is reached.  We again observe that, with the exception of the mass QoI, AMISC can approximate each QoI with the same accuracy as a single fidelity approximation at only a fraction of the cost. The AMISC approximation at some points on the convergence curves is an order of magnitude less expensive. This approaches the limit of potential acceleration, as there is only a factor of 20 difference in cost between the cheapest and most expensive models. 

\rev{From performing parameter sweeps of each QoI using our MISC approximation} we observed that mass is very close to a linear function of the random variables. Consequently, a level 1 single-fidelity sparse grid with 81 points can accurately represent mass (top left).  \rev{81 is the minimum number of points required to start the adaptive algorithm (Algorithm \ref{alg:dim-adaptivity}). The 81 points are used to compute the error indicators $\gamma_{\aindex,\bindex}^i$ for the lowest fidelity model, indexed by $\aindex=(0,0)$. These initial points are comprised of one point at the center of the domain and 2 points at the upper and lower bounds of each of the 40 random parameters.}

The AMISC approximation \rev{of mass} achieves the same final error as the single fidelity approximation, but it is slightly less accurate for moderate cost (up to 600 equivalent high-fidelity samples). However, it is more flexible than the single-fidelity approach in that AMISC can produce reasonably accurate approximations using less than 10 equivalent high-fidelity samples. This would be important if our computational budget was more restricted than the 3400 high-fidelity samples used here.

Note that the convergence curves of the approximation of thrust (top right) is noisy in the right tail. This is related to the tolerance we use for the non-linear CFD solver; in particular, the error is saturating at the error dictated by this tolerance. We could also include the solver tolerance as a physical discretization parameter, but this is left for future work.

The AMISC approximation of the thermal-layer failure criteria (lower right) is cheaper to construct for accuracies up to $10^{-4}$. For higher accuracy, the single-fidelity approximation is cheaper to construct. In the latter situation, the two coarsest structural meshes cannot be used to reduce the total error, which is smaller than the deterministic prediction error of those models. Consequently, only evaluating the highest fidelity structural model can reduce the error. Due to the assumption that the sparse grid index set $\jset$ is downward closed, the AMISC procedure requires that, for a given value of the discretization parameter $\ai_2=k$, the models parameterized by all $\ai_2=j<k$ must be evaluated before the model indexed by $k$ is evaluated. \rev{This results in redundant computation, in the form of low-fidelity evaluations} when only the model $\ai_2=k$ can be used to reduce the error of the approximation. For this same reason, the single-level sparse grid approximation of the load-layer failure criteria (bottom left) will eventually become more efficient than the AMISC approximation if we were to increase the number of samples used to generate the convergence curves.

Note that for any amount of work, the adaptive strategy we use to build both the AMISC and single-fidelity approximations attempts to balance the error in each QoI. This is because we use a worst-case error indicator which refines the approximation so that the largest error among QoI is reduced. Over the entire evolution of the adaptive algorithm, the refinement of both the single-fidelity and AMISC approximations is primarily driven by the need to reduce error in the temperature layer criteria, since the error in the approximation of this QoI is much larger than the error in the other QoI for a given cost. For example, the error of the AMISC approximation of the temperature-layer failure criteria is approximates $4\times10^{-2}$ when work is 100, whereas the error in the load-layer QoI is $2\times10^{-3}$ and the errors in the mass mass and thrust are even smaller. Consequently, at this point in the sparse grid evolution, work will be assigned to reduce the error in the load layer criteria until its error becomes comparable to the errors in the other QoI. 

\begin{figure}[htb]
\begin{center}
\includegraphics[width=\textwidth]{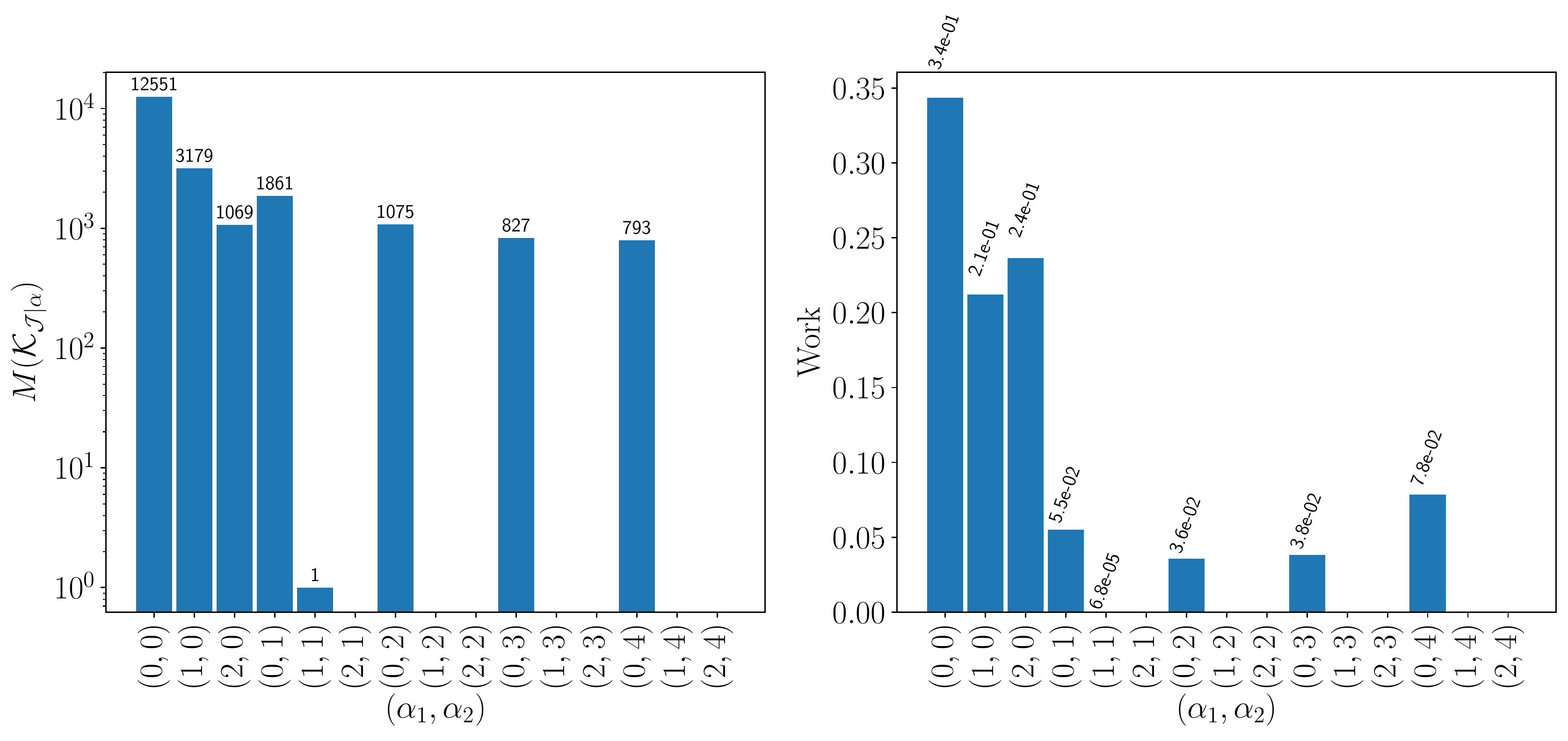}
\end{center}
\caption{The number of simulations (left) and fraction of total work (right) assigned to each nozzle model. \rev{Here $M_\kset=\mathrm{card}\left(\grid_{\kset}\right)$ is the number of samples, of the random variables, at which $\surr_\aindex$ is evaluated. The number on top of the bars in the left plot are the exact values of $M_\kset$ and the numbers on the bars on the right plot are the fraction of total work
assigned to each model indexed by $\aindex$. Finally the labels of the horizontal axis represent the indices $\aindex$ to
which AMISC has assigned at least one evaluation.}}
\label{fig:nozzle-cost-profile}
\end{figure}

In Figure \ref{fig:nozzle-cost-profile}, we plot the number of model evaluations allocated to each model discretization. The AMISC algorithm determined that the CFD and structural model are both large sources of deterministic prediction error. This is reflected by the large number of samples allocated to the most refined CFD and structural models. It is interesting to note, however, that it is not necessary to refine both the CFD and structural meshes simultaneously. There is only one sample allocated to a model, $\aindex=(1,1)$, for which both the CFD and the structural mesh have been refined . The error in the CFD model output, which is used as input to the structural model, does not significantly affect the prediction of the structural quantities of interest -- the load layer and temperature layer failure criteria. From assessing the values of $\derr^i/\dwork$ 
in the final grid, we determined that refinement of the CFD discretization parameter was driven mainly by the error in the mass and thrust. And not surprisingly, refinement of the structural discretization parameter was primarily driven by errors in the two structural QoI.

\subsubsection{Uncertainty quantification}
\label{sec-6-2-2}
In this section, we use the AMISC approximation to quantify uncertainty in the nozzle model. Specifically we report sensitivity metrics and construct PDFs for each QoI. 

\begin{figure}[htb]
\begin{center}
\includegraphics[width=\textwidth]{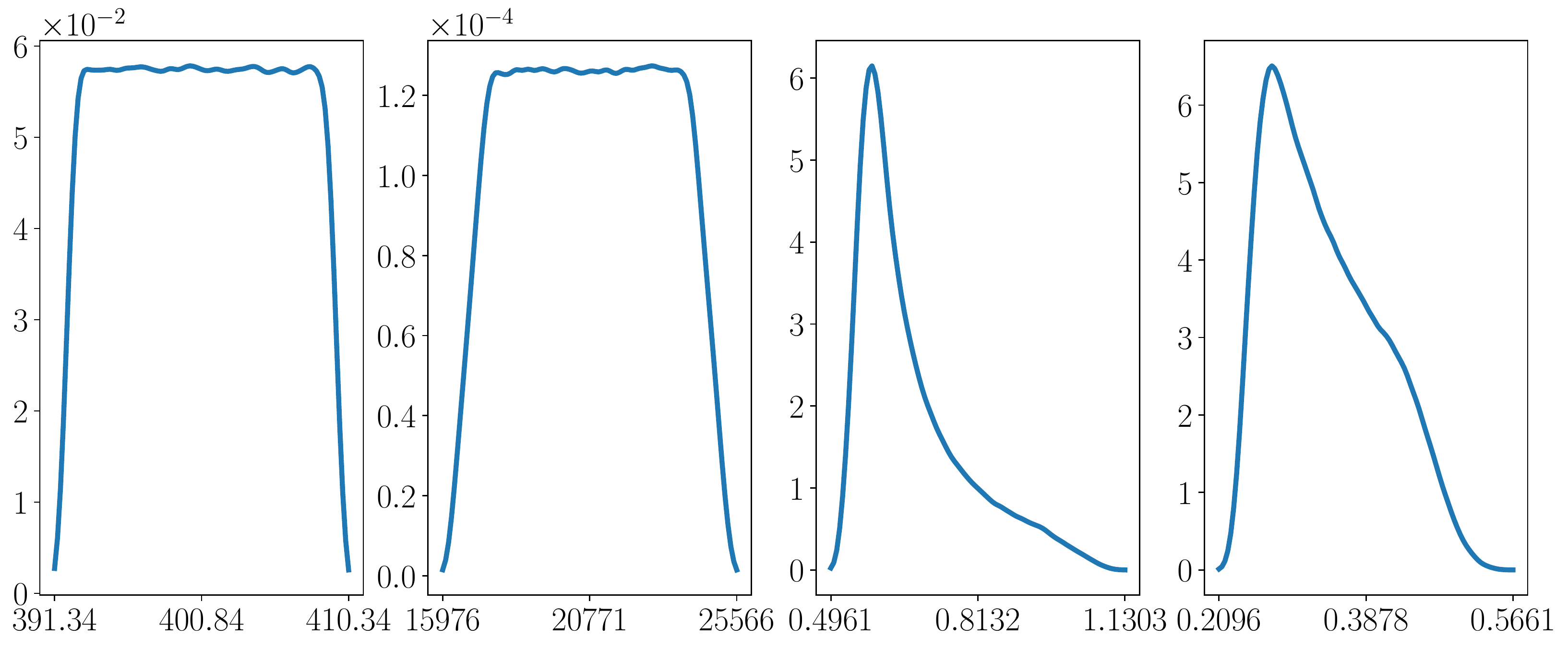}
\end{center}
\caption{PDFs for (left to right) mass, thrust, load layer failure and temperature layer failure.}
\label{fig:qoi-pdfs}
\end{figure}

In Figure \ref{fig:qoi-pdfs}, we plot the marginal PDFs of each quantity of interest. The PDFs are approximated using Gaussian kernel density estimation using $10^6$ approximate values of the QoI obtained by randomly sampling the AMISC approximation according to the distribution of the random parameters. The PDFs of both mass and thrust appear uniform. This suggests that, because the inputs are uniform, the input-output map for these QoI are primarily linear. This conclusion is also supported by plots of 1D cross sections of the map, which appear linear visually (plots omitted for brevity). The PDFs of both failure-based QoI concentrate mass on low values for each response.  Since we deem failure to occur when the failure criteria exceed 1, there is a relatively small chance of failure in the load layer and numerically zero chance of failure in the temperature layer.

\begin{figure}[htb]
\begin{center}
\includegraphics[width=\textwidth]{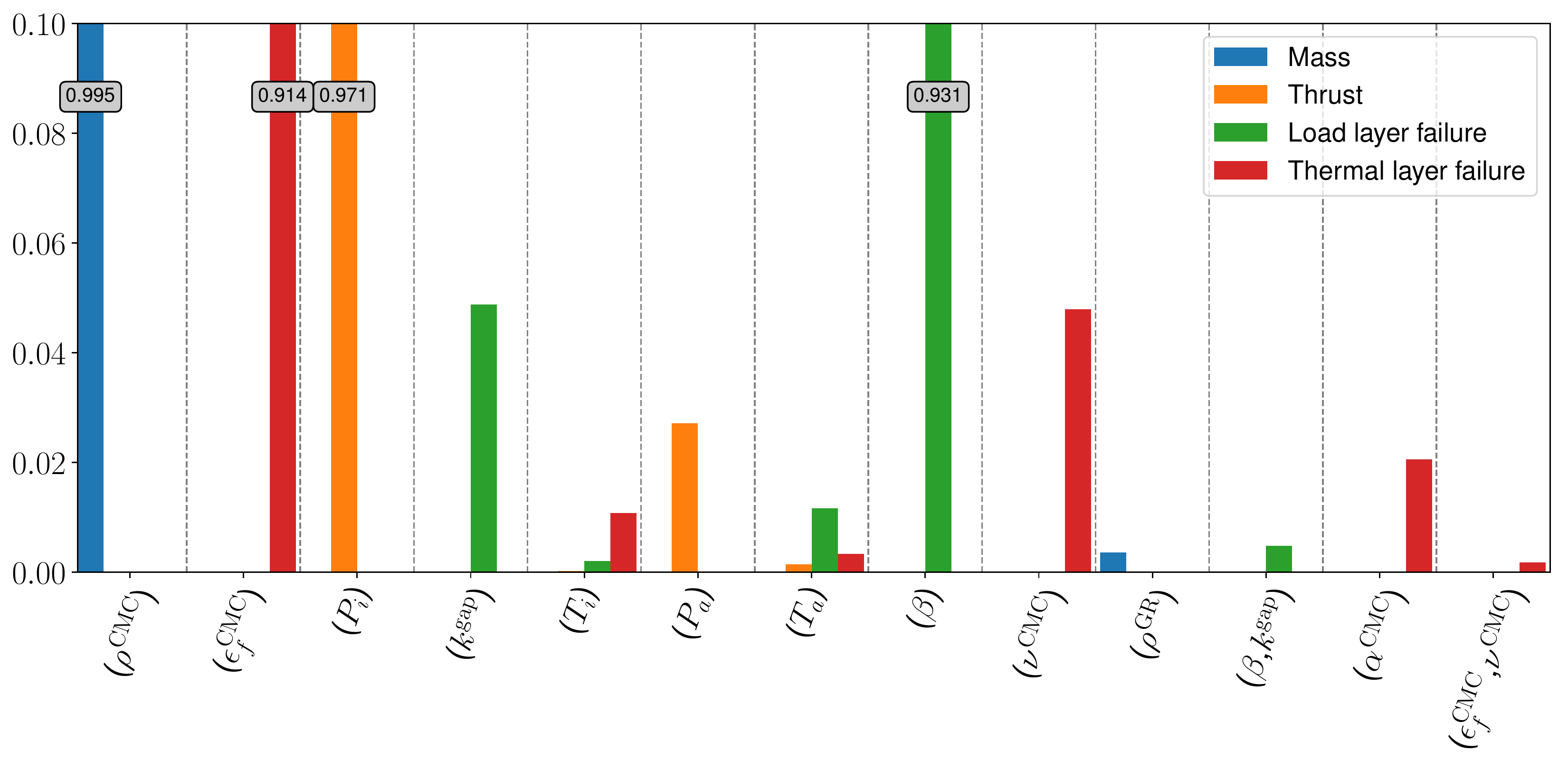}
\end{center}
\caption{Sobol indices that contribute a portion of the first 99.9\% of the variance of at least one nozzle QoI.  The labels of the x-axis represent the parameter of the Sobol index. The maximum value of any Sobol index is one. We zoom in on the y-scale to better highlight the smaller Sobol index values, and report values greater than the 0.1 plot scale in the gray boxes. The definition of the parameter symbols can be found in Appendix \ref{sec:appendix-random-variables}.}
\label{fig:sobol-indices}
\end{figure}

In Figure \ref{fig:sobol-indices}, we plot the dominant sources of uncertainty in the nozzle predictions. Specifically, we identify the parameters and parameter combinations that have large Sobol indices. The parameter combinations depicted contribute a portion of the first 99.9\% of variance of at least one nozzle QoI. It is interesting to note that, for each QoI, one parameter acting independently of all other parameters contributes over 90\% of the total variance. Moreover almost all contributions to 99.9\% of the variance of each QoI are caused by individual parameters. Only two second-order interactions are significant. However, the convergence curves \rev{in Figure~\ref{fig:nozzle-error-convergence}} suggest that resolving multivariate interactions becomes important as we drive the surrogate error to tight tolerances. This is confirmed by further investigation of high-order Sobol indices that are small but non-zero. We omit these values for brevity.

The most sensitive parameter identified for each QoI matches physical intuition. The density of the CMC heat layer $\rho^\mathrm{CMC}$ is the largest contribution to the uncertainty in the mass of the nozzle, consistent with the observations that the minimum density of this layer is an order of magnitude larger than the maximum density of the other layers and the volume of this layer is greater than any other layer. The inlet stagnation pressure $P_i$ has the greatest influence on thrust. The heat transfer coefficient $\beta$ which transfers heat from the nozzle to the environment has the greatest impact on the Load layer temperature failure ratio, which is the ratio of the temperature in the load layer to the maximum allowable temperature. When $\beta$ is decreased, less heat is able to escape the nozzle and thus the temperature in the nozzle increases. Finally, the failure strain in the CMC heat layer $\epsilon_f^\mathrm{CMC}$ has the greatest impact on the failure of the thermal failure.

The sensitivity analysis presented here suggests that uncertainty in the model can be reduced by obtaining data that constrains uncertainty on three parameters, $\rho^\mathrm{CMC}$, $\beta$ and $\epsilon_f^\mathrm{CMC}$. The inlet stagnation temperature is a function of the flight conditions and so its variance can not be constrained unless the aircraft mission is modified. Collecting experimental data and using Bayesian inference to infer updated parameter distributions to reduce uncertainty in model predictions is a topic for future work.

\section{Conclusions}
\label{sec-7}

We have described an integrated sparse grid approach for uncertainty \rev{quantification} and sensitivity analysis, which generates a \rev{surrogate model} that spans both random dimensions and model resolution dimensions.  Novel aspects of this work include the extension from moment estimation to function approximation, extension to an adaptive multi-index construction (from MISC to AMISC), and demonstration of more complete UQ workflows for realistic applications, including computation of global sensitivity analysis and probability density functions for a jet engine nozzle application.

\rev{AMISC takes a possible set of fidelity levels and greedily allocates samples to the models which provide the most predictive utility per unit cost. The set of possible models must be specified before the algorithm is initiated. Given a vector of values for each discretization parameter of a model, we defined the set of possible models to be those indexed by the Cartesian product of these vectors. As an example, the advection diffusion model considered in this paper had 3 discretization parameters which controlled the accuracy and cost of the simulation. Each one of these hyper-parameters could take 6 values resulting in $6^3=216$ possible model fidelities.

Given a Cartesian product of possible model discretizations, AMISC applies adaptive generalized sparse grid refinement to select the tensor product approximation, from among multiple grid refinement candidates, that maximizes benefit per unit cost. Reflecting a focus on complex engineering applications with multiple statistics of interest, we used an adaptive criterion based upon a convex combination of mean and variance, across a vector of QoI, to guide refinement.

In the advection-diffusion model problem presented, AMISC was able to effectively allocate resources to the downward closed subset of the $216$ model fidelities. The total computational cost required by AMISC to construct a surrogate with specified accuracy was over two orders of magnitude smaller relative to a sparse grid applied to only the most resolved model.  Moreover, an instructive transition from early resolution of stochastic error to late resolution of deterministic discretization bias is evident in the convergence history, demonstrating how multiple error sources are managed by the algorithm.}

In the nozzle example, we show a variable level of improvement in efficiency for the four different QoI, where the efficiency is strongly dependent on the nonlinearity of the response.  With greater nonlinearity, the AMISC approach displays strong savings.  However, savings are not evident for nearly linear functions due to the additional overhead of carrying multiple models; rather, benefits are limited to greater flexibility in constructing rough approximations for low cost targets. Some inefficiencies were also caused by the downward-closed index set constraint for cases where final reductions in deterministic discretization error could only be addressed using the most resolved model, resulting in some loss in relative benefits for the tails of the convergence histories.  Finally, we also demonstrate several downstream UQ products for the jet nozzle problem, including Sobol indices 
and probability density functions for each QoI, generated by processing and interrogating the final AMISC approximations.

\rev{In this paper, the relative cost of each model was determined using a small number of pilot runs prior to using AMISC. In future work, we intend to modify AMISC to determine relative model cost online during the evolution of the algorithm. Similarly, this paper focused on the case when the number of possible model discretizations is fixed a priori, as was necessary in the nozzle application. 
  That is, we consider the highest-fidelity model to be most accurate and our aim is to efficiently construct a surrogate of the statistics of this model by including evaluations of a structured set of lower-fidelity models.  In the future, we intend to extend our approach to consider an infinite set of possible discretizations, which converge to the exact solution to the governing equations. 
}
Finally, we intend to explore
approaches that alter the downward-closed requirement, in order to allow additional refinement freedom.  
This is expected to be particularly helpful when model discrepancies do not demonstrate a smooth, \rev{monotonic} decay.

\section{Acknowledgments}
\label{sec-8}
All authors were supported by the DARPA program for Enabling Quantification of Uncertainty in Physical Systems (EQUiPS). Sandia National Laboratories is a multi-mission laboratory managed and operated by National Technology and Engineering Solutions of Sandia, LLC., a wholly owned subsidiary of Honeywell International, Inc., for the U.S. Department of Energy's National Nuclear Security Administration under contract DE-NA-0003525. The views expressed in the article do not necessarily represent the views of the U.S. Department of Energy or the United States Government.

\ifx \journal\undefined
\bibliographystyle{plain} 
\fi

\bibliography{references}

\clearpage \appendix

\section{Nozzle model input uncertainties}
\label{sec-9}
\label{sec:appendix-random-variables}
The materials of the nozzle, inlet conditions, atmospheric conditions, and heat transfer are all sources of uncertainty that effect predictions of our nozzle model. In this section we explicitly identify and characterize (assign a probability distribution) to each source of uncertainty.

\subsection{Material Constants}
\label{sec-9-1}
\begin{wrapfigure}{r}{0.4\textwidth}
    \centering
    \includegraphics[width=0.5\linewidth]{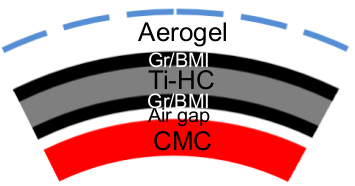}
        \caption{Individual layers in the nozzle wall.}
   \label{fig:components}
\end{wrapfigure}

Figure \ref{fig:components} shows the individual wall layers in the nozzle wall. The material properties of each of these layers \rev{a}ffect nozzle performance.  The thermal insulation ``aerogel'' is not modeled, but all other layers are modeled.
An incomplete set of material properties is gathered from a variety of experimental data sources. 

\noindent Missing macroscopic properties are either determined through fundamental physical analysis of the material's structure or are approximated by values for similar materials. 
Table \ref{tab:matPropCMC} provides the properties used in the \texttt{MULTI-F} analysis for the CMC material, Table \ref{tab:matPropGR-BMI} for the GR-BMI material, and Table \ref{tab:matPropTI-HC} for the TI-HC material.
Lastly, Table \ref{tab:airgapmatPropPanel} lists the thermal conductivity of the air gap and gives the distribution for the assumed yield stress of the composite panel material used for the baffles.

\floatsetup[table]{capposition=top}

\begin{table}[hb]
\caption[CMC material properties]{Assumed \textit{isotropic} material properties for the heat layer's ceramic matrix composite material (CMC).}
\label{tab:matPropCMC}
\begin{center}
\begin{tabular}[]{ c| c | c H | c }
\textbf{Property} & \textbf{symbol} & \textbf{Units} & \textbf{Nominal Value} & \textbf{Distribution} \\ \hline
Density & $\rho^\mathrm{CMC}$ & $\frac{\textrm{kg}}{\textrm{m}^3}$ & 2410 & $\mathcal{U}(2293.789,2496.486)$ \\ \hline
Elastic modulus & $E^\mathrm{CMC}$ & GPa & 67.1 & $\mathcal{U}(58.940,76.163)$ \\ \hline
Poisson ratio & $\nu^\mathrm{CMC}$ & & 0.33 & $\mathcal{U}(0.23,0.43)$\\ \hline
Thermal conductivity & $k^\mathrm{CMC}$ & $\frac{\textrm{W}}{\textrm{m-K}}$ & 1.41 & $\mathcal{U}(1.37,1.45)$ \\ \hline
Thermal expansion coefficient & $\alpha^\mathrm{CMC}$ & $\textrm{K}^{-1} \times 10^{-6}$ & 0.24 & $\mathcal{U}(0.228,0.252)$ \\ \hline \hline
Max service temperature & $T_{\max}^{\mathrm{CMC}}$ & K & 973 & $\mathcal{U}(963,983)$ \\ \hline
Failure strain & $\epsilon_{f}^\mathrm{CMC}$ & \% & 0.07 & $\mathcal{U}(0.049788,0.096441)$ \\ \hline
\end{tabular}
\end{center}
\end{table}

\begin{table}[htb]
\caption[Gr-BMI material properties]{Macroscopic laminate material properties for the load layer's graphite/bismaleimide (GR-BMI) composite layers used in the axisymmetric shell specification.}
\label{tab:matPropGR-BMI}
\begin{center}
\begin{tabular}[]{ c | c | c H | c }
\textbf{Property} & \textbf{Symbol} & \textbf{Units} & \textbf{Nominal Value} & \textbf{Distribution} \\ \hline
Density & $\rho^\mathrm{GR}$ & $\frac{\textrm{kg}}{\textrm{m}^3}$ & 1568 & $\mathcal{U}(1563,1573)$ \\ \hline
Elastic moduli & $E_{1} = E_{2}$ & GPa & 60 & $\mathcal{U}(57,63)$ \\ \hline
In-plane shear modulus & $G_{12}$ & GPa & 23.31 & $\mathcal{U}(22.6,24.0)$ \\ \hline
Poisson ratios & $\nu_{12} = \nu_{21}$ & & 0.344 & $\mathcal{U}(0.334,0.354)$\\ \hline
Mutual influence coef (first kind) & $\mu_{1,12} = \mu_{2,12}$ & & 0.0 & $\mathcal{U}(-0.1,0.1)$ \\ \hline
Mutual influence coef (second kind) & $\mu_{12,1} = \mu_{12,2}$ & & 0.0 & $\mathcal{U}(-0.1,0.1)$ \\ \hline
Thermal conductivity & $k_{1} = k_{2}$ & $\frac{\textrm{W}}{\textrm{m-K}}$ & 3.377 & $\mathcal{U}(3.208,3.546)$ \\ \hline
Thermal conductivity & $k_{3}$ & $\frac{\textrm{W}}{\textrm{m-K}}$ & 3.414 & $\mathcal{U}(3.243,3.585)$ \\ \hline
Thermal expansion coef & $\alpha_{1} = \alpha_{2}$ & $\textrm{K}^{-1} \times 10^{-6}$ & 1.200 & $\mathcal{U}(1.16,1.24)$ \\ \hline
Thermal expansion coef & $\alpha_{12}$ & $\textrm{K}^{-1} \times 10^{-6}$ & 0.0 & $\mathcal{U}(-0.04,0.04)$ \\ \hline \hline
Max service temperature & $T_{\max}^{\mathrm{GR}}$ & K & 505 & $\mathcal{U}(500,510)$ \\ \hline
Failure strain (tension) & $\epsilon_{f,1}^t = \epsilon_{f,2}^t$ & \% & 0.75 & $\mathcal{U}(0.675,0.825)$ \\ \hline
Failure strain (compression) & $\epsilon_{f,1}^c = \epsilon_{f,2}^c$ & \% & -0.52 & $\mathcal{U}(-0.572,-0.494)$ \\ \hline
Failure strain (shear) & $\gamma_{f}$ & \% & 0.17 & $\mathcal{U}(0.153,0.187)$ \\ \hline
\end{tabular} \\
\end{center}
\end{table}

\begin{table}[htb]
\caption[TI-HC material properties]{Assumed \textit{isotropic} macroscopic material properties for titanium honeycomb layer (TI-HC).}
\label{tab:matPropTI-HC}
\begin{center}
\begin{tabular}[]{ c | c | c H | c }
\textbf{Property} & \textbf{Symbol} & \textbf{Units} & \textbf{Nominal Value} & \textbf{Distribution} \\ \hline
Density & $\rho^\mathrm{TI}$ & $\frac{\textrm{kg}}{\textrm{m}^3}$ & 179.57 & $\mathcal{U}(177.77,181.37)$ \\ \hline
Elastic modulus& $E^\mathrm{TI}$ & GPa & 1.90 & $\mathcal{U}(1.587,2.823)$ \\ \hline
Poisson ratio & $\nu^\mathrm{TI}$ & & 0.178 & $\mathcal{U}(0.160,0.196)$ \\ \hline
Thermal conductivity & $k^\mathrm{TI}$ & $\frac{\textrm{W}}{\textrm{m-K}}$ & 0.708 &  $\mathcal{U}(0.680,0.736)$ \\ \hline
Thermal expansion coefficient & $\alpha^\mathrm{TI}$ & $\textrm{K}^{-1} \times 10^{-6}$ & 2.97 & $\mathcal{U}(2.88,3.06)$ \\ \hline \hline
Max service temperature & $T_{\max}^{\mathrm{TI}}$ & K & 755 &  $\mathcal{U}(745,765)$ \\ \hline
Yield stress & $\sigma_{Y}^{\mathrm{TI}}$ & MPa & 12.9 & $\mathcal{U}(9.676,16.951)$ \\ \hline
\end{tabular}
\end{center}
\end{table}

\begin{table}[htb]
\caption[Air gap material and panel structure material properties]{Assumed properties of air gap between thermal and load layers and of the panel structure used in baffles.}
\label{tab:airgapmatPropPanel}
\begin{center}
\begin{tabular}[]{ c | c | c H | c }
\textbf{Property} & \textbf{Symbol} & \textbf{Units} & \textbf{Nominal Value} & \textbf{Distribution} \\ \hline
Thermal conductivity & $k^\mathrm{gap}$ & $\frac{\textrm{W}}{\textrm{m-K}}$ & 0.0425 &  $\mathcal{U}(0.0320,0.0530)$\\ \hline
Yield strength & $\sigma_{Y,B}$ & MPa & 324 & $\mathcal{U}(56.876,99.219)$ \\ \hline
\end{tabular}
\end{center}
\end{table}

\subsection{Mission Parameters}
\label{sec-9-2}

We focus on a typical reconnaissance mission for a small high-subsonic unmanned military aircraft. The mission includes climbing at maximum rate to a cruise altitude of 43,000 ft, cruising at Mach 0.92 for a specified distance to an observation point (e.g., 500 km), loitering at an altitude of 43,000 ft and Mach 0.5 for 2 hours, and then returning to the takeoff point. On the return, the aircraft descends to 10,000 ft and cruises at Mach 0.9 in a high-speed ``dash'' segment lasting several kilometers before landing.
The analysis of this mission showed that the climb segment was the most critical for nozzle performance since maximum thrust was required at all altitudes, leading to the highest temperatures and pressures at the inlet of the nozzle. In particular, the state of climb right before beginning the cruise segment was the most critical in terms of stresses and temperatures experienced by the nozzle. Consequently, we focus on the climb segment of the mission here. Table \ref{tab:missionParams} summarizes the mission parameters used in the nozzle analysis. We assume the altitude and mach number are fixed at 40,000 ft and 0.511, respectively.

\begin{table}[htb]
\caption[Mission parameters]{Mission parameters.}
\label{tab:missionParams}
\begin{center}
\begin{tabular}[]{ c |c | c | H c }
\textbf{Parameter} & \textbf{Symbol} & \textbf{Units} & \textbf{Nominal Value} & \textbf{Distribution} \\ \hline
Inlet stagnation pressure & $P_i$ &Pa & 97,585 & $\mathcal{U}(86362.122,113062.433)$ \\ \hline
Inlet stagnation temperature & $T_i$ & K & 955.0 & $\mathcal{U}(928.729,981.600)$ \\ \hline
Atmospheric pressure & $P_a$ & Pa & 18,754 & $\mathcal{U}(17386.668,20206.125)$ \\ \hline
Atmospheric temperature & $T_a$ & K & 216.7 & $\mathcal{U}(202.856,231.297)$ \\ \hline
Heat transfer coefficient to environment & $\beta$ & $\frac{\textrm{W}}{\textrm{m}^2\textrm{-K}}$ & 12.62 & $\mathcal{U}(4.043,37.376)$ \\ \hline
\end{tabular}
\end{center}
\end{table}
\end{document}

%% file: defs.tex
\newcommand{\rev}[1]{#1}

\newcommand{\V}[1]{{\boldsymbol{#1}}}
\newcommand{\C}[1]{{\mathcal{#1}}}
\newcommand{\var}[1]{{\mathbb{V}[{#1}]}}
\newcommand{\mean}[1]{{\mathbb{E}[{#1}]}}

\def\nv{d}  
\def\na{{n_\alpha}}  
\def\ai{\alpha}    
\def\aindex{{\V{\ai}}}  
\def\nb{{n_\beta}}  
\def\bi{\beta}    
\def\bindex{{\V{\bi}}}  
\def\dv{x}
\def\dvv{{\V{\dv}}}    
\def\rv{z}         
\def\rvv{{\V{\rv}}}  
\def\rV{Z}         
\def\rVv{{\V{\rV}}}  
\def\rvvsupp{\Gamma} 
\def\nq{\rev{q}}     
\def\reals{\mathbb{R}} 
\def\indexset{\mathcal{I}}
\def\grid{\mathcal{Z}}                
\def\values{\mathcal{F}}              
\newcommand{\samp}[1]{\rvv^{(#1)}} 
\newcommand{\unisamp}[2]{\rv_{#1}^{(#2)}} 
\def\surr{\hat{f}} 
\def\sgrid{\surr_{\aindex,\iset}} 
\newcommand{\isgrid}[1]{\surr_{\aindex,\iset(#1)}} 
\def\ffem{\surr_\aindex}
\newcommand{\imisc}[1]{\surr_{\jset(#1)}}
\def\misc{\surr_{\jset}}

\def\physdom{D}

\def\derr{\Delta E_{\aindex,\bindex}}
\def\dwork{\Delta W_{\aindex,\bindex}}
\newcommand{\norm}[2]{\lVert #1 \rVert_{#2}}
\newcommand{\mynorm}[1]{\norm{#1}{L^p(\rvvsupp)}}
\newcommand{\abs}[1]{\lvert #1 \rvert}
\newcommand{\dx}[1]{\;\mathrm{d}#1}
\def\basis{\phi}

\def\diffgrid{\grid^{\mathrm{diff}}}

\def\iset{\mathcal{I}}
\def\jset{\mathcal{J}}
\def\kset{{\mathcal{K}_{\jset\mid\aindex}}}
\def\cL{\jset}
\def\cA{\mathcal{A}}
\def\bl{{\aindex,\bindex}}
\def\blb{{[\bl]}}
\def\bls{{\aindex^\star,\bindex^\star}}
\def\blsb{{[\bls]}}

\DeclareMathOperator*{\argmax}{argmax}

\newcommand{\myint}[1]{\int_\rvvsupp #1 \dx{\omega(\rvv)}}

\def\pbwt{\omega}

\newcommand*{\bigchi}{\mbox{\Large$\chi$}}